\documentclass[12pt,sort&compress]{elsarticle}
\usepackage{tikz, amsmath,url}
\usetikzlibrary{decorations.markings}
\usepackage{listings}
\lstdefinelanguage{maple}
	{morekeywords={true, false, try, catch, return, break, error, 
	               module, export, local, option, in, use,
                 and, or, not, xor, xnor,
                 if, then, elif, else, fi,
                 while, for, from, by, to, do, od,
                 proc, nargs, local, global, end, NULL}}
\newdefinition{definition}{Definition}
\newtheorem{theorem}{Theorem}
\newtheorem{corollary}{Corollary}
\newdefinition{example}{Example}
\newproof{proof}{Proof}

\newcommand{\cross}{\operatorname{cr}}
\newcommand{\nest}{\operatorname{ne}}
\newcommand{\NCN}{\operatorname{NCN}}
\newcommand{\NN}{\operatorname{NN}}
\newcommand{\NC}{\operatorname{NC}}
\newcommand{\arc}{\operatorname{Arc}}

\newcommand{\figurefontsize}{\small}
\newcommand{\tableaufontsize}{\relax}
\providecommand{\texorpdfstring}[2]{#1} 
\tikzstyle{pnt}=[circle,fill,inner sep=1.5pt]
\tikzstyle{state}=[draw, circle,inner sep=1pt,minimum size=6mm]
\tikzstyle{uni}=[decoration={markings,%
   mark=at position 0.3 with {\arrow[line width=1pt]{>};}},%
   postaction={decorate}]
\tikzstyle{unir}=[decoration={markings,%
   mark=at position 0.3 with {\arrowreversed[line width=1pt]{>};}},%
   postaction={decorate}]
\tikzstyle{loop above right}=[above right,out=60,in=30,loop]
\tikzstyle{loop above left}=[above left,out=150,in=120,loop]
\tikzstyle{loop below left}=[below right,out=240,in=210,loop]
\tikzstyle{loop below right}=[below right,out=330,in=300,loop]

\begin{document}

\begin{frontmatter}
\title{Crossings and Nestings  for Arc-Coloured Permutations}
\author{Lily Yen}
\address{Dept.\ of Math.\ \& Stats., Capilano University, North Vancouver, B.C., Canada;\\{\normalfont also:} Dept.\ of Math., Simon Fraser University, Burnaby, B.C., Canada}

\begin{abstract}
The equidistribution of many crossing and nesting statistics exists in several combinatorial objects like matchings, set partitions, permutations, and embedded labelled graphs. The involutions switching nesting and crossing numbers for set partitions given by Krattenthaler, also by Chen, Deng, Du, Stanley, and Yan, and for permutations given by Burrill, Mishna, and Post involved passing through tableau-like objects. Recently, Chen and Guo for matchings, and Marberg  for set partitions extended the result to coloured arc annotated diagrams. We prove that symmetric joint distribution continues to hold for arc-coloured permutations. As in Marberg's recent work, but through a different interpretation, we also conclude that the ordinary generating functions for all $j$-noncrossing, $k$-nonnesting, $r$-coloured permutations according to size $n$ are rational functions.  We use the interpretation to automate  the generation of these rational series for both noncrossing and nonnesting coloured set partitions and permutations.
\end{abstract}
\end{frontmatter}\thispagestyle{empty}

\section{Introduction}

Crossing and nesting statistics have intrigued combinatorialists for many decades. For example, it is well known that Catalan numbers, 
\(
c_n = \frac{1}{n+1} \binom{2n}{n}
\), count the number of noncrossing matchings on $[2n]$ which is also  the number of nonnesting matchings of the same size. The concept of crossing and nesting was then extended to higher numbers where symmetric joint distribution continues to hold not only for matchings~\citep{GB89}, but also for set partitions~\citep{Chetal07, Kratt06}, labelled graphs~\citep{deMi07}, set partitions of classical types~\citep{RuSt10}, and permutations~\citep{BuMiPo10}. In all cases, bijective proofs were given; and for some, generating functions were found. 

Inspired by recent works of \citet{ChG11} on coloured matchings and \citet{Mar12} on coloured set partitions, we give a bijection to establish symmetric joint distribution of crossing and nesting statistics for \emph{arc-coloured permutations}. We also  show that the ordinary generating functions for $j$-noncrossing, $k$-nonnesting, $r$-coloured permutations according to size $n$ are  rational functions.

\subsection{Definitions and Terminology}
A permutation $S$ of the set $[n]:=\{ 1, 2, \dots, n\}$ is a bijection from $[n]$ to itself,  $\sigma: [n] \rightarrow [n]$. Using two-line notation, 
\(
   S= \left(\begin{smallmatrix}
   1               & 2             & 3                   & \dots & n   \\
   \sigma(1)  &\sigma(2)  & \sigma(3)        & \dots  &\sigma(n)
   \end{smallmatrix} \right)
\).
An \emph{arc annotated diagram} is a labelled graph on $n$ vertices drawn horizontally, labelled left to right consecutively such that $\arc(i,j)$ joins vertex $i$ to vertex $j$. A permutation has a representation as an arc annotated diagram where $\arc(i, \sigma(i))$ is drawn as an upper arc for $\sigma(i) \ge i$, and a lower arc for $\sigma(i) < i$. Note that the dissymmetry draws a fixed point in $S$ as an upper loop. When this diagram is restricted to only the upper arcs (or lower arcs) with all $n$ vertices, then it also represents a set partition of $[n]$. Separately, we call these\emph{ upper} and \emph{ lower arc diagrams} of a permutation.  From this diagram, we define a \emph{$k$-crossing} (resp.\ \emph{$k$-nesting}) as $k$ arcs 
\(
\{(i_1, j_1), (i_2, j_2), \dots, (i_k, j_k)\}
\) all mutually cross, or 
\(
i_1 < i_2 < \dots < i_k < j_1 < j_2 < \dots < j_k
\) (resp.\ nest, i.\,e.\ 
\(
i_1 < i_2 < \dots < i_k < j_k < j_{k-1}< \dots < j_1
\)
) as shown in Figure~\ref{fig:kcrossing} (resp.\ Figure~\ref{fig:knesting}). We also need a variant: \emph{enhanced $k$-crossing} (resp.\ \emph{enhanced $k$-nesting}) where 
\(
i_1 < i_2 < \dots < i_k \le j_1 < j_2 < \dots < j_k
\)
 (resp.\ 
 \(
 i_1 < i_2 < \dots < i_k \le j_k < j_{k-1}< \dots < j_1
 \)
 ) as shown in Figure~\ref{fig:kenhancedcr} (resp.\ Figure~\ref{fig:kenhancedne}).

\begin{figure}\centering
\begin{minipage}[b]{0.45\linewidth}
\figurefontsize\centering
\begin{tikzpicture}[bend left=35,scale=0.8]
   \node[pnt,label=below:$i_1$] at (1,0) {};
   \node[pnt,label=below:$i_2$] at (2,0) {};
   \node[label=below:$\dots$] at (3,0) {};
   \node[pnt,label=below:$i_k$] at (4,0) {};
   \node[pnt,label=below:$j_1$] at (5,0) {};
   \node[pnt,label=below:$j_2$] at (6,0) {};
   \node[label=below:$\dots$] at (7,0) {};
   \node[pnt,label=below:$j_k$] at (8,0) {};
   \draw (1,0) to (5,0);
   \draw (2,0) to (6,0);
   \draw (3,0) to (7,0);
   \draw (4,0) to (8,0);	
\end{tikzpicture}
\caption{The arc diagram of a $k$-crossing}
\label{fig:kcrossing}
\end{minipage}
\quad
\begin{minipage}[b]{0.45\linewidth}
\figurefontsize\centering
\begin{tikzpicture}[bend left=30,scale=0.8]
   \node[pnt,label=below:$i_1$] at (1,0) {};
   \node[pnt,label=below:$i_2$] at (2,0) {};
   \node[label=below:$\dots$] at (3,0) {};
   \node[pnt,label=below:$i_k$] at (4,0) {};
   \node[pnt,label=below:$j_k$] at (5,0) {};
   \node[label=below:$\dots$] at (6,0) {};
   \node[pnt,label=below:$j_2$] at (7,0) {};
   \node[pnt,label=below:$j_1$] at (8,0) {};
   \draw (1,0) to (8,0);
   \draw (2,0) to (7,0);
   \draw (3,0) to (6,0);
   \draw (4,0) to (5,0);
\end{tikzpicture}
\caption{The arc diagram of a $k$-nesting}
\label{fig:knesting}
\end{minipage}
\end{figure}

\begin{figure}
\centering
\begin{minipage}[b]{0.45\linewidth}
\figurefontsize\centering
\begin{tikzpicture}[scale=0.8,bend left=35]
   \node[pnt,label=below:$i_1$] at (1,0) {};
   \node[pnt,label=below:$i_2$] at (2,0) {};
   \node[label=below:$\dots$] at (3,0) {};
   \node[pnt,label={below:$i_k=j_1$}] at (4,0) {};
   \node[pnt,label=below:$j_2$] at (5,0) {};
   \node[label=below:$\dots$] at (6,0) {};
   \node[pnt,label=below:$j_k$] at (7,0) {};
   \draw (1,0) to (4,0);
   \draw (2,0) to (5,0);
   \draw (3,0) to (6,0);
   \draw (4,0) to (7,0);	
\end{tikzpicture}
\caption{The arc diagram of an \emph{enhanced} $k$-crossing}
\label{fig:kenhancedcr}
\end{minipage}
\quad
\begin{minipage}[b]{0.45\linewidth}
\figurefontsize\centering
\begin{tikzpicture}[scale=0.8,bend left=30]
   \node[pnt,label=below:$i_1$] at (1,0) {};
   \node[pnt,label=below:$i_2$] at (2,0) {};
   \node[label=below:$\dots$] at (3,0) {};
   \node[pnt,label={below:$i_k=j_k$}] at (4,0) {};
   \node[label=below:$\dots$] at (5,0) {};
   \node[pnt,label=below:$j_2$] at (6,0) {};
   \node[pnt,label=below:$j_1$] at (7,0) {};
   \draw (1,0) to (7,0);
   \draw (2,0) to (6,0);
   \draw[bend left=45] (3,0) to (5,0);
   \draw[loop right] (4,0) to ();
\end{tikzpicture}
\caption{The arc diagram of an \emph{enhanced} $k$-nesting}
\label{fig:kenhancedne}
\end{minipage}
\end{figure}

We need both notions of crossings and nestings for permutations because the \emph{enhanced} definitions are used for upper arc diagrams whereas the other definitions (without \emph{enhanced}), for lower arc diagrams. This is in accordance with the literature~\citep{Corteel07} on permutation statistics for weak exceedances and pattern avoidance. We define the \emph{crossing number}, $\cross(S) = j$ (resp.\  \emph{nesting number}, $\nest(S) = k$) of a permutation $S$ as the maximum $j$ (resp.\ $k$) such that $S$ has a $j$-enhanced crossing (resp.\ $k$-enhanced nesting) in the upper arc diagram or a $j$-crossing (resp.\ $k$-nesting) in the lower arc diagram. When a permutation $S$ does not have a $j$-(enhanced)-crossing (resp.\ $k$-(enhanced)-nesting), then we say $S$ is $j$-noncrossing (resp.\ $k$-nonnesting). \citet*{BuMiPo10} gave an involution mapping between the set of  permutations of $[n]$ with $\cross(S)=j$ and $\nest(S)=k$ and  those with $\cross(S)=k$ and $\nest(S)=j$, thus extending the result of symmetric joint distribution for matchings and set partitions of \citet*{Chetal07} and \citet{Kratt06} to permutations.

Next, \citet{ChG11} generalized symmetric equidistribution of crossing and nesting statistics to \emph{coloured} complete matchings. Most recently, \citet{Mar12} extended the result to coloured set partitions with a novel way of proving that the ordinary generating functions of $j$-noncrossing, $k$-nonnesting, $r$-coloured partitions according to size $n$ are rational functions. We extend their results to $r$-arc-coloured permutations, or $r$-coloured permutations in short.

Coloured permutations are generalizations of permutations represented as arc annotated diagrams. Once the arcs are coloured to satisfy $j$-noncrossing and $k$-nonnesting conditions for each colour class, the resulting arc diagrams can be represented in the topological graph theoretic book embedding setting~\citep{Mm06}, each colour on a separate page while the vertices are on the spine of the book. The differences are two fold: each page satisfies the crossing/nesting conditions instead of finding a minimum number of pages, noncrossing on each page, to represent a given (non-planar) graph, and the number of pages is not necessarily minimal with respect to the crossing/nesting conditions. Secondary RNA structures with different bonding  energies have been analysed in the book embedding setting, naturally represented as coloured set partitions~\citep{clote2012page}; however, arc-coloured permutations have yet to find a natural application.

Some caution on terminology is in order here. Group properties of col\-our\-ed permutations have been widely studied since the $1990$'s~\citep{Borodin99, Stump12}, but there the colours are assigned to \emph{vertices} instead of \emph{arcs}.

\subsection{Main Theorem}

  Since crossing and nesting statistics involves arcs, we define an \emph{$r$-coloured permutation} parallel to~\citep{Mar12} as a pair, $(S, \phi)$ consisting of a permutation of $[n]$ and an arc-colour assigning map $\phi: \arc(S) \rightarrow [r]$, and use a capital Greek letter, $\Sigma$, to denote these objects. We say $\Sigma$ has a $k$-crossing (resp.\ $k$-nesting) if $k$ arcs of the \emph{same} colour cross (resp.\ nest). As always throughout this paper, \emph{enhanced} statistics is applied to \emph{upper} arc diagrams while \emph{non-enhanced} for \emph{lower} arc diagrams of permutations. Define $\cross(\Sigma)$ (resp.\ $\nest(\Sigma)$) as the maximum integer $k$ such that $\Sigma$ has a $k$-crossing (resp.\ $k$-nesting). The bijection of~\citep{BuMiPo10} can be extended to establish symmetric joint distribution of the numbers $\cross(\Sigma)$ and $\nest(\Sigma)$ over $r$-coloured permutations preserving opener and closer sequences (equivalently, sets of minimal and maximal elements of each block when upper arc and lower arc diagrams are viewed separately as set partitions). 
  
  More formally, vertices of a permutation are of five types,
   an opener (
\begin{tikzpicture}
      \draw[bend left] (1.2,0) node[pnt]{} to +(0.3,0.2);
      \draw[bend right] (1.2, 0) to (1.5, -0.2);
\end{tikzpicture}
),
 a closer (
  \begin{tikzpicture}
      \draw[bend right] (1.5,0) node[pnt]{} to (1.2, 0.2);
      \draw[bend left] (1.5, 0) to (1.2, -0.2);
   \end{tikzpicture}
 ), 
 a fixed point (
 \begin{tikzpicture}
      \draw[loop above] (0,0) node[pnt]{} to (0.0);
   \end{tikzpicture}
   ) , 
   an upper transitory (
   \begin{tikzpicture}
      \draw[bend left] (0,0) node[pnt]{} to +(0.3,0.2);
      \draw[bend right] (0,0) to (-0.3, 0.2);
   \end{tikzpicture}
   ), and 
   a lower transitory(
   \begin{tikzpicture}
      \draw[bend right] (0,0) node[pnt]{} to (0.3,-0.2);
      \draw[bend left] (0,0) to (-0.3, -0.2);
   \end{tikzpicture}). For a particular $\Sigma$, restricting to only one colour, both upper arc and lower arc diagrams can be seen as set partitions whose minimal block elements are the openers, and maximal block elements are the closers. For upper arc diagrams, both a fixed point and an upper transitory contribute to  the set of minimal (opener) and the set of maximal (closer) elements over blocks of the set partition. Lower arc diagrams are set partitions in Marberg's partition setting, thus Theorem1.1 and Corollary 1.2 of~\citep{Mar12} apply exactly here.
  
  Given an $r$-coloured permutation $\Sigma = (S, \phi)$, let the set of openers (resp.\ the set of closers) be $\mathcal{O}(\Sigma)$ (resp.\ $\mathcal{C}(\Sigma)$) of the uncoloured permutation, $S$. For all positive integers, $j$ and $k$, and subsets $O$, $C \subseteq [n]$, define 
  \(
  \NCN_{j,k}^{O,C} (n,r)
  \)
   to be the number of $r$-coloured permutations $\Sigma$ of $[n]$ with $\cross(\Sigma) < j$, $\nest(\Sigma) < k$, $\mathcal{O}(\Sigma) = O$, and $\mathcal{C}(\Sigma) = C$. Then Theorem~\ref{thm:jointjk} is analogous to Theorem 1.1 in~\citep{Chetal07, Mar12} for $r$-coloured permutations.
  
 \begin{theorem}
 \label{thm:jointjk}
 For all positive integers, $j$ and $k$, and subsets $O$, $C \subseteq [n]$, 
  \(
 \NCN_{j,k}^{O,C} (n,r) = \NCN_{k,j}^{O,C} (n,r)
 \).
 \end{theorem}
 As customary in the literature, we let $\NCN_{j,k}(n,r)$ denote the number of all $r$-coloured, $j$-noncrossing, $k$-nonnesting permutations of $[n]$. Summing both sides of Theorem~\ref{thm:jointjk} over all $O$, $C \subseteq [n]$ gives the generalization of~\citep{Chetal07, Mar12} for Corollary~\ref{cor:NC}. We also let $\NC_k(n,r)$ (resp. $\NN_k(n,r)$) denote the number of $k$-noncrossing (resp. $k$-nonnesting) $r$-coloured permutations on $[n]$.

\begin{corollary}
\label{cor:NC}
For all integers, $j, k, n, r$, 
\(
\NCN_{j,k}(n,r)  = \NCN_{k,j}(n,r)
\)
 and $\NC_k(n,r) = \NN_k(n,r)$.
\end{corollary}

\subsection{Plan}

The tools needed for the proof of Theorem~\ref{thm:jointjk} are given in Section~\ref{sec:background}.  Section~\ref{sec:pfMT} gives the proof of Theorem~\ref{thm:jointjk} combining essential ingredients  of both~\citep{BuMiPo10, Mar12} with the added care of managing both upper and lower arc diagrams simultaneously where both notions of crossing and nesting are applied. The transfer matrix approach Marberg used to establish the rationality of the ordinary generating function, 
\(
\sum_{n \ge 0} \NCN_{j,k}(n+1,r) x^n
\)
 for set partitions of size $n+1$ is through translating the original problem to counting all closed walks of $n$-steps with certain column and row length restrictions (according to $j, k$) for each component from 
 \(
 \mathbf{\emptyset} \in \mathbf{Y}^r
 \), 
 that is, $r$ copies of the Hasse diagram of the Young lattice. This idea cannot be extended to permutations on $(\mathbf{Y}^r, \mathbf{Y}^r)$ because upper arc diagrams are dependent on lower arc diagrams. However, another interpretation of Marberg's multigraphs $\mathcal{G}_{j,k,r}$ in terms of the types of vertices and colours of edges leads to the multigraphs for $r$-coloured permutations which permits the application of transfer matrix method to draw the same conclusion: The ordinary generating function, 
 \(
 \sum_{n \ge 0} \NCN_{j,k}(n,r) x^n
 \)
  for $j$-noncrossing, $k$-nonnesting, $r$-coloured permutations is rational. The combination of the method of generating trees and finite state automata in the interpretation can be extended to other  combinatorial objects where both crossing and nesting statistics are bounded, thus leading to the same conclusion that the corresponding generating functions are rational.

\section{Background}
\label{sec:background}
The proof of Theorem~\ref{thm:jointjk} requires working knowledge of the theory of integer partition, especially its representation as Young diagrams, the Hasse diagram of the Young lattice, and the RSK-algorithm for filling  positive integers to obtain  the beginning of some standard Young tableau. We refer the reader to Volume $2$ of Stanley's Enumerative Combinatorics~\citep{Stanley99} for more details.

Define a partition of $n\in \mathbf{N}$ to be a sequence 
\(
\lambda = (\lambda_1, \lambda_2, \dots, \lambda_k) \in \mathbf{N}^k
\)
 such that 
\(
\sum_{i=1}^k \lambda_i = n
\), and 
\(
\lambda_1 \ge \lambda_2 \ge \dots \ge \lambda_k
\).
 If $\lambda$ is a partition of $n$, we write $\lambda \vdash n$ or $| \lambda | = n$. The non-zero terms $\lambda_i$ are called the parts of $\lambda$, and we say $\lambda$ has $k$ parts if $\lambda_k > 0$. We can draw $\lambda$ using a left-justified array of boxes with $\lambda_i$ boxes in row $i$. For example, 
 \(
 \lambda = (5, 3, 2, 2, 1)
 \)
  is drawn as 
	\begin{tikzpicture}[scale=0.1]
		\draw (0,4) grid (5,5);
		\draw (0,3) grid (3,4);
		\draw (0,1) grid (2,3);
		\draw (0,0) grid (1,1);
	\end{tikzpicture}.
This representation is the Young diagram of a partition. To ``add a box'' to a partition $\lambda$ means to obtain a partition $\mu$ such that $| \lambda| + 1 = | \mu|$, and $\lambda$'s Young diagram is included in that of $\mu$. This inclusion induces a partial order on the set of partitions of non-negative integers, denoted by $\mathbf{Y}$, or the Young lattice. When we place integers $1, 2, \dots, n$ in all $n$ boxes of a Young diagram so that entries increase in each row and column, we produce a standard Young tableau, abbreviated as SYT. As one builds an SYT from the empty set through the process of adding a box at a time, a sequence of integer partitions, 
\(
( \lambda^0 = \emptyset, \lambda^1, \lambda^2, \dots, \lambda^n)
\)
 emerges where $\lambda^{i-1} \subset \lambda^i$, 
 and $|\lambda^i| = |\lambda^{i-1}| + 1$. In addition to adding a box, we include ``deleting a box'' and ``doing nothing'' for the following four types in Definition~\ref{dfn:T}.

\begin{definition}
\label{dfn:T}
We define four types of sequences of tableaux,  
\(
T = ( \lambda^0 = \emptyset, \lambda^1, \lambda^2, \dots, \lambda^n)
\), 
where $\lambda^0 = \lambda^n = \emptyset$ such that $\lambda^i$ is obtained from $\lambda^{i-1}$ for each $i \in [n]$ by one of the three actions: adding a box, deleting a box, or doing nothing.
\begin{enumerate}
\item A \emph{semi-oscillating tableau} is any such sequence $T$. 
\item An \emph{oscillating tableau} has distinct neighbouring $\lambda^i$'s.
\item A \emph{vacillating tableau} is any such sequence $T$ which has  $\lambda^{i-1} \subseteq \lambda^i$ when $i$ is even, and  $\lambda^{i-1} \supseteq \lambda^i$ when $i$ is odd.
\item  A \emph{hesitating tableau} is any such sequence $T$ which has  $\lambda^{i-1} \subseteq \lambda^i$ when $i$ is odd, and  $\lambda^{i-1} \supseteq \lambda^i$ when $i$ is even.
\end{enumerate}
\end{definition}

In the uncoloured case, \citet{Mar12} links the sequence $T$ to an $n$-step walk on the Hasse diagram of the Young lattice, $\mathbf{Y}$ where ``doing nothing'' is also counted as a step. 
For his enumeration purposes,  Marberg's definitions differ slightly from~\citep{Chetal07} to achieve that these $n$-step walks are closed walks from $\emptyset$. Though we will not walk on an ordered pair of $r$-tuple Hasse diagrams, we will keep the requirement that each sequence $T$ begins and ends with $\emptyset$.

\section{Proof of Main Theorem}
\label{sec:pfMT}
To warm-up for the proof of Theorem~\ref{thm:jointjk}, we give  examples of how three of the four types of $T$'s from Definition~\ref{dfn:T} are used to encode different combinatorial objects. Two local rules for changing set partitions to involutions are needed for the last two tableaux in Section~\ref{ex:T}:  Rule H for hesitating tableaux tracking \emph{enhanced} statistics in upper arcs and Rule V for vacillating tableaux. 
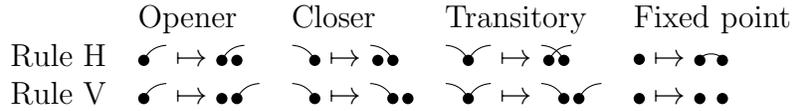
\begin{center}
\begin{tabular}{lllll}
& Opener& Closer& Transitory& Fixed point\\

Rule H& 
   \begin{tikzpicture}
      \draw[bend left] (0,0) node[pnt]{} to +(0.3,0.2);
   \end{tikzpicture}%
   ${}\mapsto{}$%
   \begin{tikzpicture}
      \draw (1.5,0) node[pnt]{};
      \draw[bend left] (1.3,0) node[pnt]{} to +(0.3,0.2);
   \end{tikzpicture}
&
   \begin{tikzpicture}
      \draw[bend right] (0,0) node[pnt]{} to (-0.3,0.2);
   \end{tikzpicture}%
   ${}\mapsto{}$%
   \begin{tikzpicture}
      \draw (1.3,0) node[pnt]{};
      \draw[bend right] (1.5,0) node[pnt]{} to (1.2,0.2);
   \end{tikzpicture}
&
   \begin{tikzpicture}
      \draw[bend left] (0,0) node[pnt]{} to +(0.3,0.2);
      \draw[bend right] (0,0) to (-0.3, 0.2);
   \end{tikzpicture}%
   ${}\mapsto{}$%
   \begin{tikzpicture}
      \draw[bend right] (1.5,0) node[pnt]{} to (1.2, 0.2);
      \draw[bend left] (1.3,0) node[pnt]{} to +(0.3,0.2);
   \end{tikzpicture}
&
   \begin{tikzpicture}
      \draw[loop above] (0,0) node[pnt]{};
   \end{tikzpicture}%
   ${}\mapsto{}$%
   \begin{tikzpicture}
      \draw[bend left = 80] (1,0) to (1.3, 0);
      \draw (1,0) node[pnt]{};
      \draw (1.3, 0) node[pnt]{};
    \end{tikzpicture}
\\
Rule V& 
   \begin{tikzpicture}
      \draw[bend left] (0,0) node[pnt]{} to +(0.3,0.2);
   \end{tikzpicture}%
   ${}\mapsto{}$%
   \begin{tikzpicture}
      \draw (1.3,0) node[pnt]{};
      \draw[bend left] (1.5,0) node[pnt]{} to +(0.3,0.2);
   \end{tikzpicture}
&
   \begin{tikzpicture}
      \draw[bend right] (0,0) node[pnt]{} to (-0.3,0.2);
   \end{tikzpicture}%
   ${}\mapsto{}$%
   \begin{tikzpicture}
      \draw (1.5,0) node[pnt]{};
      \draw[bend right] (1.3,0) node[pnt]{} to (1,0.2);
   \end{tikzpicture}
&
   \begin{tikzpicture}
      \draw[bend left] (0,0) node[pnt]{} to +(0.3,0.2);
      \draw[bend right] (0,0) to (-0.3, 0.2);
   \end{tikzpicture}%
   ${}\mapsto{}$%
   \begin{tikzpicture}
      \draw[bend right] (1.3,0) node[pnt]{} to (1, 0.2);
      \draw[bend left] (1.5,0) node[pnt]{} to +(0.3,0.2);
   \end{tikzpicture}
&
   \begin{tikzpicture}
      \draw (0,0) node[pnt]{};
   \end{tikzpicture}%
   ${}\mapsto{}$%
   \begin{tikzpicture}
      \draw (1,0) node[pnt]{};
      \draw (1.3,0) node[pnt]{};
   \end{tikzpicture}
\end{tabular}
\end{center}

\subsection{Examples of three  encodings}
\label{ex:T}
\begin{example}
\label{ex:soin} A \emph{semi-oscillating tableau sequence} encoding an involution.
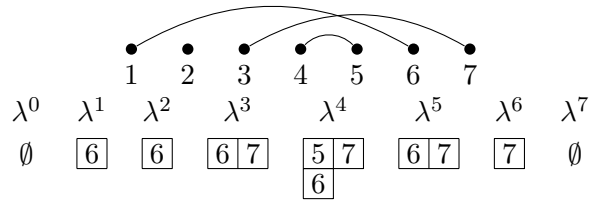
\begin{center}\figurefontsize
	\begin{tikzpicture}[scale=0.75]
		 \foreach \i in {1,...,7}
       	 \node[pnt,label=below:$\i$] at (\i,0)(\i) {};
       \draw(1)  to [bend left=30] (6);
      \draw(3)  to [bend left=30] (7);
      \draw(4)  to [bend left=45] (5);
   \end{tikzpicture}
\\
\begin{tabular}{cccccccc}
$\lambda^0$ & $\lambda^1$ & $\lambda^2$ & $\lambda^3$ & $\lambda^4$ & $\lambda^5$  & $\lambda^6$ & $\lambda^7$  \\
      \begin{tikzpicture}[scale=0.4,baseline={(0,0.4)}]
            \node at (0.5,1.5) {$\emptyset$};
      \end{tikzpicture} &
	\begin{tikzpicture}[scale=0.4,baseline={(0,0.4)}]
		\draw (0,1) grid (1,2);
		\node at (0.5,1.5) {\tableaufontsize$6$};
	\end{tikzpicture} &
	 \begin{tikzpicture}[scale=0.4,baseline={(0,0.4)}]
		\draw (0,1) grid (1,2);
		\node at (0.5,1.5) {\tableaufontsize$6$};
	\end{tikzpicture} &
	\begin{tikzpicture}[scale=0.4,baseline={(0,0.4)}]
		\draw (0,1) grid (2,2);
		\node at (0.5,1.5) {\tableaufontsize$6$};
		\node at (1.5,1.5) {\tableaufontsize$7$};
	\end{tikzpicture} &
	\begin{tikzpicture}[scale=0.4,baseline={(0,0.4)}]
		\draw (0,1) grid (2,2);
		\draw (0,0) grid (1,1);
		\node at (0.5,1.5) {\tableaufontsize$5$};
		\node at (0.5,0.5) {\tableaufontsize$6$};
		\node at (1.5,1.5) {\tableaufontsize$7$};
	\end{tikzpicture} &
	\begin{tikzpicture}[scale=0.4,baseline={(0,0.4)}]
		\draw (0,1) grid (2,2);
		\node at (0.5,1.5) {\tableaufontsize$6$};
		\node at (1.5,1.5) {\tableaufontsize$7$};
	\end{tikzpicture} &
	\begin{tikzpicture}[scale=0.4,baseline={(0,0.4)}]
		\draw (0,1) grid (1,2);
		\node at (0.5,1.5) {\tableaufontsize$7$};
	\end{tikzpicture} &
      \begin{tikzpicture}[scale=0.4,baseline={(0,0.4)}]
            \node at (0.5,1.5) {$\emptyset$};
      \end{tikzpicture} 
\end{tabular}
\end{center}
\end{example}

\begin{example}
\label{ex:vsp} A \emph{vacillating tableau sequence} encoding a set partition.
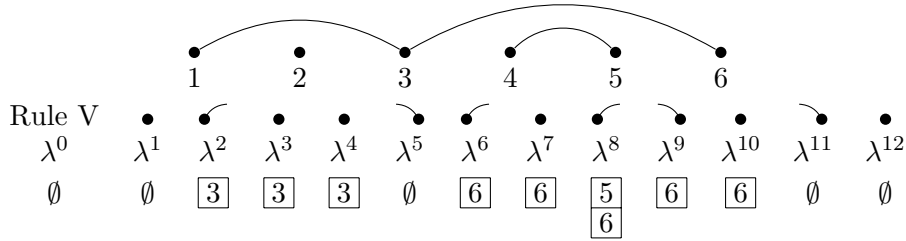
\begin{center}
\figurefontsize
	\begin{tikzpicture}[scale=1.4]
		 \foreach \i in {1,...,6}
       	 \node[pnt,label=below:$\i$] at (\i,0)(\i) {};
       \draw(1)  to [bend left=30] (3);
      \draw(3)  to [bend left=30] (6);
      \draw(4)  to [bend left=45] (5);
   \end{tikzpicture}
\\
\begin{tabular}{*{13}{c}}
Rule V &
   \begin{tikzpicture}
      \draw (1,0) node[pnt]{};
   \end{tikzpicture}
   & 
   \begin{tikzpicture}
      \draw[bend left] (1.2,0) node[pnt]{} to +(0.3,0.2);
   \end{tikzpicture}
   & 
    \begin{tikzpicture}
      \draw (1,0) node[pnt]{};
   \end{tikzpicture}
   & 
    \begin{tikzpicture}
      \draw (1.2,0) node[pnt]{};
   \end{tikzpicture}
   & 
   \begin{tikzpicture}
      \draw[bend right] (1.3,0) node[pnt]{} to (1, 0.2);
   \end{tikzpicture}
   & 
   \begin{tikzpicture}
      \draw[bend left] (1.5,0) node[pnt]{} to +(0.3,0.2);
   \end{tikzpicture}
   & 
   \begin{tikzpicture}
      \draw (1,0) node[pnt]{};
   \end{tikzpicture}
   & 
   \begin{tikzpicture}
      \draw[bend left] (1.2,0) node[pnt]{} to +(0.3,0.2);
   \end{tikzpicture}
   & 
  \begin{tikzpicture}
      \draw[bend right] (1.3,0) node[pnt]{} to (1,0.2);
   \end{tikzpicture}
   &
   \begin{tikzpicture}
      \draw (1.5,0) node[pnt]{};
   \end{tikzpicture}
   &
   \begin{tikzpicture}
      \draw[bend right] (1.3,0) node[pnt]{} to (1,0.2);
   \end{tikzpicture}
   &
   \begin{tikzpicture}
      \draw (1.5,0) node[pnt]{};
   \end{tikzpicture}
       \\
    $\lambda^0$ & $\lambda^1$ & $\lambda^2$ & $\lambda^3$ & $\lambda^4$ & $\lambda^5$  & $\lambda^6$ & $\lambda^7$   & $\lambda^8$  & $\lambda^9$  & $\lambda^{10}$   & $\lambda^{11}$  & $\lambda^{12}$
  \\
	\begin{tikzpicture}[scale=0.4,baseline={(0,0.4)}]
            \node at (0.5,1.5) {$\emptyset$};
      \end{tikzpicture} &
	\begin{tikzpicture}[scale=0.4,baseline={(0,0.4)}]
            \node at (0.5,1.5) {$\emptyset$};
      \end{tikzpicture} &
	\begin{tikzpicture}[scale=0.4,baseline={(0,0.4)}]
		\draw (0,1) grid (1,2);
		\node at (0.5,1.5) {\tableaufontsize$3$};
	\end{tikzpicture} &
	\begin{tikzpicture}[scale=0.4,baseline={(0,0.4)}]
		\draw (0,1) grid (1,2);
		\node at (0.5,1.5) {\tableaufontsize$3$};
	\end{tikzpicture} &
	\begin{tikzpicture}[scale=0.4,baseline={(0,0.4)}]
		\draw (0,1) grid (1,2);
		\node at (0.5,1.5) {\tableaufontsize$3$};
	\end{tikzpicture} &
	\begin{tikzpicture}[scale=0.4,baseline={(0,0.4)}]
            \node at (0.5,1.5) {$\emptyset$};
      \end{tikzpicture} &
	\begin{tikzpicture}[scale=0.4,baseline={(0,0.4)}]
		\draw (0,1) grid (1,2);
		\node at (0.5,1.5) {\tableaufontsize$6$};
	\end{tikzpicture} &
	\begin{tikzpicture}[scale=0.4,baseline={(0,0.4)}]
		\draw (0,1) grid (1,2);
		\node at (0.5,1.5) {\tableaufontsize$6$};
	\end{tikzpicture} &
	\begin{tikzpicture}[scale=0.4,baseline={(0,0.4)}]
		\draw (0,0) grid (1,2);
		\node at (0.5,1.5) {\tableaufontsize$5$};
		\node at (0.5,0.5) {\tableaufontsize$6$};
	\end{tikzpicture} &
	\begin{tikzpicture}[scale=0.4,baseline={(0,0.4)}]
		\draw (0,1) grid (1,2);
		\node at (0.5,1.5) {\tableaufontsize$6$};
	\end{tikzpicture} &
	\begin{tikzpicture}[scale=0.4,baseline={(0,0.4)}]
		\draw (0,1) grid (1,2);
		\node at (0.5,1.5) {\tableaufontsize$6$};
	\end{tikzpicture} &
	\begin{tikzpicture}[scale=0.4,baseline={(0,0.4)}]
            \node at (0.5,1.5) {$\emptyset$};
      \end{tikzpicture} &
	\begin{tikzpicture}[scale=0.4,baseline={(0,0.4)}]
            \node at (0.5,1.5) {$\emptyset$};
      \end{tikzpicture} 
\end{tabular}
\end{center}
\end{example}

\begin{example}
\label{ex:hesp} A \emph{hesitating tableau sequence} encoding an \emph{enhanced} set partition.
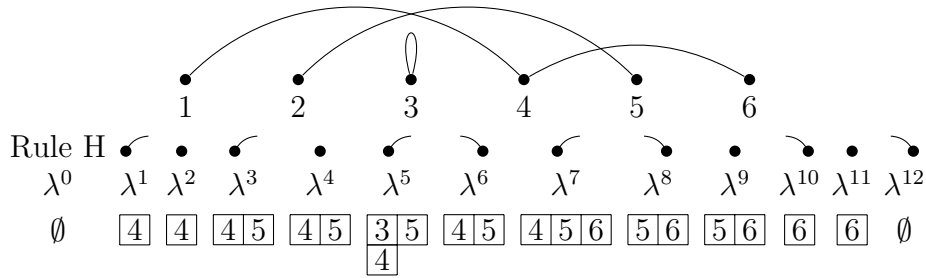
\begin{center}
	\begin{tikzpicture}[scale=1.5]
		 \foreach \i in {1,...,6}
       	 \node[pnt,label=below:$\i$] at (\i,0)(\i) {};
       \draw(1)  [bend left=45] to (4);
      \draw(4)  to [bend left=30] (6);
      \draw(2)  to [bend left=45] (5);
      \draw[loop above] (3) to ();
   \end{tikzpicture}
\\
\setlength{\tabcolsep}{2pt}%
\begin{tabular}{*{13}{c}}
Rule H &
   \begin{tikzpicture}
      \draw[bend left] (1,0) node[pnt]{} to +(0.3,0.2);
   \end{tikzpicture}
   & 
   \begin{tikzpicture}
      \draw (1.2,0) node[pnt]{};
   \end{tikzpicture}
   &
   \begin{tikzpicture}
      \draw[bend left] (1,0) node[pnt]{} to +(0.3,0.2);
   \end{tikzpicture}
   & 
   \begin{tikzpicture}
      \draw (1.2,0) node[pnt]{};
   \end{tikzpicture}
   &
   \begin{tikzpicture}
       \draw[bend left] (1,0) node[pnt]{} to +(0.3, 0.2);
    \end{tikzpicture}   
   & 
   \begin{tikzpicture}
       \draw[bend right] (1.5,0) node[pnt]{} to (1.2, 0.2);
    \end{tikzpicture}   
   & 
    \begin{tikzpicture}
      \draw[bend left] (1.3,0) node[pnt]{} to +(0.3,0.2);
   \end{tikzpicture}
    & 
   \begin{tikzpicture}
      \draw[bend right] (1.5,0) node[pnt]{} to (1.2, 0.2);
   \end{tikzpicture}
    & 
     \begin{tikzpicture}
      \draw (1.3,0) node[pnt]{};
   \end{tikzpicture}
   &
   \begin{tikzpicture}
      \draw[bend right] (1.5,0) node[pnt]{} to (1.2,0.2);
   \end{tikzpicture}
   &
   \begin{tikzpicture}
      \draw (1.3,0) node[pnt]{};
   \end{tikzpicture}
   &
   \begin{tikzpicture}
      \draw[bend right] (1.5,0) node[pnt]{} to (1.2,0.2);
   \end{tikzpicture}
\\
    $\lambda^0$ & $\lambda^1$ & $\lambda^2$ & $\lambda^3$ & $\lambda^4$ & $\lambda^5$  & $\lambda^6$ & $\lambda^7$   & $\lambda^8$  & $\lambda^9$  & $\lambda^{10}$   & $\lambda^{11}$  & $\lambda^{12}$
\\
	\begin{tikzpicture}[scale=0.4,baseline={(0,0.4)}]
            \node at (0.5,1.5) {$\emptyset$};
      \end{tikzpicture} &
	\begin{tikzpicture}[scale=0.4,baseline={(0,0.4)}]
		\draw (0,1) grid (1,2);
            \node at (0.5,1.5) {\tableaufontsize$4$};
	\end{tikzpicture} &
	\begin{tikzpicture}[scale=0.4,baseline={(0,0.4)}]
		\draw (0,1) grid (1,2);
            \node at (0.5,1.5) {\tableaufontsize$4$};
	\end{tikzpicture} &
	\begin{tikzpicture}[scale=0.4,baseline={(0,0)}]
		\draw (0,0) grid (2,1);
            \node at (0.5,0.5) {\tableaufontsize$4$};
            \node at (1.5,0.5) {\tableaufontsize$5$};
	\end{tikzpicture} &
	\begin{tikzpicture}[scale=0.4,baseline={(0,0)}]
		\draw (0,0) grid (2,1);
            \node at (0.5,0.5) {\tableaufontsize$4$};
            \node at (1.5,0.5) {\tableaufontsize$5$};
	\end{tikzpicture} &
	\begin{tikzpicture}[scale=0.4,baseline={(0,0.4)}]
		\draw (0,0) grid (1,1);
		\draw (0,1) grid (2,2);
            \node at (0.5,1.5) {\tableaufontsize$3$};
            \node at (1.5,1.5) {\tableaufontsize$5$};
            \node at (0.5,0.5) {\tableaufontsize$4$};
	\end{tikzpicture} &
	\begin{tikzpicture}[scale=0.4,baseline={(0,0)}]
		\draw (0,0) grid (2,1);
            \node at (0.5,0.5) {\tableaufontsize$4$};
            \node at (1.5,0.5) {\tableaufontsize$5$};
	\end{tikzpicture} &
	\begin{tikzpicture}[scale=0.4,baseline={(0,0)}]
		\draw (0,0) grid (3,1);
            \node at (0.5,0.5) {\tableaufontsize$4$};
            \node at (1.5,0.5) {\tableaufontsize$5$};
            \node at (2.5,0.5) {\tableaufontsize$6$};
	\end{tikzpicture} &
	\begin{tikzpicture}[scale=0.4,baseline={(0,0)}]
		\draw (0,0) grid (2,1);
            \node at (0.5,0.5) {\tableaufontsize$5$};
            \node at (1.5,0.5) {\tableaufontsize$6$};
	\end{tikzpicture} &
	\begin{tikzpicture}[scale=0.4,baseline={(0,0)}]
		\draw (0,0) grid (2,1);
            \node at (0.5,0.5) {\tableaufontsize$5$};
            \node at (1.5,0.5) {\tableaufontsize$6$};
	\end{tikzpicture} &
	\begin{tikzpicture}[scale=0.4,baseline={(0,0.4)}]
		\draw (0,1) grid (1,2);
            \node at (0.5,1.5) {\tableaufontsize$6$};
	\end{tikzpicture} &
	\begin{tikzpicture}[scale=0.4,baseline={(0,0.4)}]
		\draw (0,1) grid (1,2);
            \node at (0.5,1.5) {\tableaufontsize$6$};
	\end{tikzpicture}&
	\begin{tikzpicture}[scale=0.4,baseline={(0,0.4)}]
            \node at (0.5,1.5) {$\emptyset$};
      \end{tikzpicture} 
\end{tabular}
\end{center}
\end{example}
Remark: All three $T$'s are constructed from vertex $1$ to the right using closer labels as fillings of SYT whereas in~\citep{Chetal07}, they are constructed from the right with opener labels as fillings of SYT. In example~\ref{ex:soin}, an opener corresponds to a new box in $T$, and its matching closer is the filling of the box. The RSK algorithm is applied when a new box (an opener) with its new filling (its matching closer) is added to the previous tableau. The tableau loses a box of a certain filling when  a matching closer is encountered in the arc diagram. The remaining labels rearrange themselves to remain an SYT whose shape is included in the previous tableau. Apply the encoding similarly for examples~\ref{ex:vsp} and~\ref{ex:hesp}. When a node \begin{tikzpicture} 
 	\draw (0,1) node[pnt]{}; 
 \end{tikzpicture}
 is present in Rule V or Rule H, that is, ``doing nothing'', no change occurs to the tableau.

\subsection{Proof of Theorem~\ref{thm:jointjk}}
\begin{proof}
We show an involution between the set of $r$-coloured permutations of $[n]$ with maximal crossing number $j$, nesting number $k$ and those with maximal crossing number $k$ and nesting number $j$.

Given an $r$-coloured permutation of $[n]$, say $\Sigma = (S, \phi)$, first consider its corresponding uncoloured permutation $S$. Let $O$ be $\mathcal{O}(S)$, the set of openers and $C$ be $\mathcal{C}(S)$, the set of closers. Applying the involution of~\citep{BuMiPo10} results in another permutation with the same $O$ and $C$ while switching maximal crossing and nesting numbers.

Now for each colour class, the resulting arc diagram is no longer a permutation, but two set partitions: \emph{enhanced} for the upper arc diagram, and  \emph{non-enhanced} for the lower arc diagram.
We employ the same encoding techniques from~\citep{BuMiPo10}:
\begin{enumerate}
\item[Step 1] Translate the upper arc diagram into a hesitating tableau sequence, and the lower arc diagram into a vacillating tableau sequence.
\item[Step 2] Perform a component-wise transpose to each tableau sequence.
\item[Step 3] Apply reverse RSK to fill each tableau in the sequence from the right to the left. 
\item[Step 4] Translate the newly filled sequence of tableaux back to arc diagrams according to its own rule.
\end{enumerate}
Thus, we obtain the resulting arc diagram with its upper and lower arc components where maximal crossing and nesting numbers are switched because the bijections of~\citep{Chetal07, Kratt06, BuMiPo10} interchange maximal column length with maximal row length  while preserving sets of maximal and minimal block elements. This interchange achieved through taking the conjugate (transpose) of each tableau translates to the switching of maximal nesting and crossing numbers while preserving the sets of openers and closers. The preservation of these sets when restricted to one colour of arcs permits the involution to be applied separately to all arcs of the same colour, one colour at a time, without interfering with the sets of openers and closers from other colour classes. Finally, the combination of all $r$ involutions, one for each colour, produces the desired $r$-coloured permutation such that for each colour, crossing number and nesting number are switched. If the original $r$-coloured $\Sigma$ is $j$-noncrossing and $k$-nonnesting, then its image after the $r$-fold involution is $j$-nonnesting and $k$-noncrossing.
\end{proof}

\subsection{An example of a \texorpdfstring{$2$}{2}-coloured permutation}
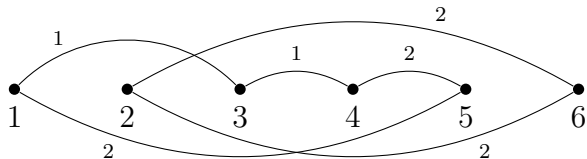
\begin{figure}[htp]
   \begin{center}
	\begin{tikzpicture}[scale=1.5]
		 \foreach \i in {1,...,6}
       	 \node[pnt,label=below:$\i$] at (\i,0)(\i) {};
       \draw[bend left=45](1)  to node[pos=0.2,above]{$\scriptstyle1$}  (3);
      \draw[bend left=30](2) to node[pos=0.7, above]{$\scriptstyle2$}  (6);
      \draw[bend left=30](4) to node[above]{$\scriptstyle2$} (5);
      \draw[bend left=30](3)  to node[above]{$\scriptstyle1$}  (4);
      \draw[bend right=30](1) to node[pos=0.2, below]{$\scriptstyle2$}  (5);
      \draw[bend right=30](2) to node[pos=0.8, below]{$\scriptstyle2$} (6);
         \end{tikzpicture}
         \caption[The image of the $2$-coloured permutation]
         {The image of  Example~\ref{ex:perm} under the involution in the proof of Theorem~\ref{thm:jointjk}}
          \label{fig:inperm}
   \end{center}
\end{figure}

We show a $2$-coloured permutation where we apply the involution of the proof of Theorem~\ref{thm:jointjk} to find its image.
\begin{example}
\label{ex:perm} A permutation encoded by a \emph{hesitating tableau sequence}, $\lambda_1$ for colour $1$, $\lambda_2$ for colour $2$ in the \emph{upper} arcs and a \emph{vacillating tableau sequence}, $\mu_2$ for colour $2$ in the \emph{lower} arcs.
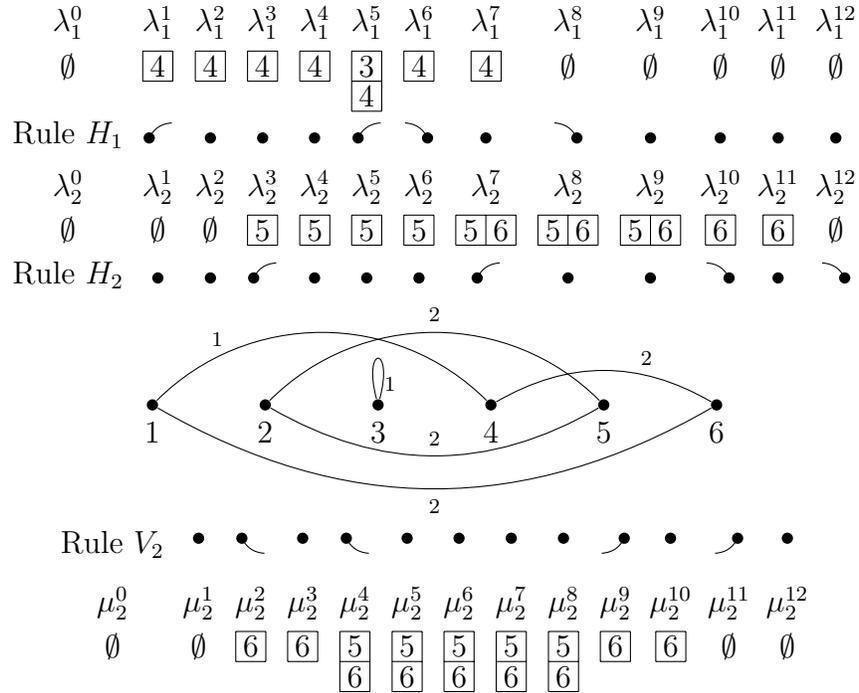
\begin{center}
\setlength{\tabcolsep}{3pt}%
\begin{tabular}{*{13}{c}}
   $\lambda_1^0$ & $\lambda_1^1$ & $\lambda_1^2$ & $\lambda_1^3$ & $\lambda_1^4$ & $\lambda_1^5$  & $\lambda_1^6$ & $\lambda_1^7$   & $\lambda_1^8$  & $\lambda_1^9$  & $\lambda_1^{10}$   & $\lambda_1^{11}$  & $\lambda_1^{12}$
\\
	\begin{tikzpicture}[scale=0.4,baseline={(0,0.4)}]
            \node at (0.5,1.5) {$\emptyset$};
      \end{tikzpicture} &
	\begin{tikzpicture}[scale=0.4,baseline={(0,0.4)}]
		\draw (0,1) grid (1,2);
            \node at (0.5,1.5) {\tableaufontsize$4$};
	\end{tikzpicture} &
	\begin{tikzpicture}[scale=0.4,baseline={(0,0.4)}]
		\draw (0,1) grid (1,2);
            \node at (0.5,1.5) {\tableaufontsize$4$};
	\end{tikzpicture} &
	\begin{tikzpicture}[scale=0.4,baseline={(0,0.4)}]
		\draw (0,1) grid (1,2);
            \node at (0.5,1.5) {\tableaufontsize$4$};
	\end{tikzpicture} &
	\begin{tikzpicture}[scale=0.4,baseline={(0,0.4)}]
		\draw (0,1) grid (1,2);
            \node at (0.5,1.5) {\tableaufontsize$4$};
	\end{tikzpicture} &	
	\begin{tikzpicture}[scale=0.4,baseline={(0,0.4)}]
		\draw (0,0) grid (1,1);
		\draw (0,1) grid (1,2);
            \node at (0.5,1.5) {\tableaufontsize$3$};
            \node at (0.5,0.5) {\tableaufontsize$4$};
	\end{tikzpicture} &
	\begin{tikzpicture}[scale=0.4,baseline={(0,0.4)}]
		\draw (0,1) grid (1,2);
            \node at (0.5,1.5) {\tableaufontsize$4$};
	\end{tikzpicture} &
	\begin{tikzpicture}[scale=0.4,baseline={(0,0.4)}]
		\draw (0,1) grid (1,2);
            \node at (0.5,1.5) {\tableaufontsize$4$};
	\end{tikzpicture} &
	\begin{tikzpicture}[scale=0.4,baseline={(0,0.4)}]
            \node at (0.5,1.5) {$\emptyset$};
      \end{tikzpicture} &
	\begin{tikzpicture}[scale=0.4,baseline={(0,0.4)}]
            \node at (0.5,1.5) {$\emptyset$};
      \end{tikzpicture}&
      \begin{tikzpicture}[scale=0.4,baseline={(0,0.4)}]
            \node at (0.5,1.5) {$\emptyset$};
      \end{tikzpicture}&
      \begin{tikzpicture}[scale=0.4,baseline={(0,0.4)}]
            \node at (0.5,1.5) {$\emptyset$};
      \end{tikzpicture}&
	\begin{tikzpicture}[scale=0.4,baseline={(0,0.4)}]
            \node at (0.5,1.5) {$\emptyset$};
      \end{tikzpicture} 
\\
Rule $H_1$ &
   \begin{tikzpicture}
      \draw[bend left] (1,0) node[pnt]{} to +(0.3,0.2);
   \end{tikzpicture}
   & 
   \begin{tikzpicture}
      \draw (1.2,0) node[pnt]{};
   \end{tikzpicture}
   &
   \begin{tikzpicture}
      \draw (1.2,0) node[pnt]{};
   \end{tikzpicture}
   &   
   \begin{tikzpicture}
      \draw (1.2,0) node[pnt]{};
   \end{tikzpicture}
   &
   \begin{tikzpicture}
       \draw[bend left] (1,0) node[pnt]{} to +(0.3, 0.2);
    \end{tikzpicture}   
   & 
   \begin{tikzpicture}
       \draw[bend right] (1.5,0) node[pnt]{} to (1.2, 0.2);
    \end{tikzpicture}   
   & 
    \begin{tikzpicture}
      \draw (1.2,0) node[pnt]{};
   \end{tikzpicture}
   &
   \begin{tikzpicture}
      \draw[bend right] (1.5,0) node[pnt]{} to (1.2, 0.2);
   \end{tikzpicture}
    & 
     \begin{tikzpicture}
      \draw (1.2,0) node[pnt]{};
   \end{tikzpicture}
   &
  \begin{tikzpicture}
      \draw (1.2,0) node[pnt]{};
   \end{tikzpicture}
   &
   \begin{tikzpicture}
      \draw (1.2,0) node[pnt]{};
   \end{tikzpicture}
   &
   \begin{tikzpicture}
      \draw (1.2,0) node[pnt]{};
   \end{tikzpicture}
\\[1ex]
   $\lambda_2^0$ & $\lambda_2^1$ & $\lambda_2^2$ & $\lambda_2^3$ & $\lambda_2^4$ & $\lambda_2^5$  & $\lambda_2^6$ & $\lambda_2^7$   & $\lambda_2^8$  & $\lambda_2^9$  & $\lambda_2^{10}$   & $\lambda_2^{11}$  & $\lambda_2^{12}$
\\
	\begin{tikzpicture}[scale=0.4,baseline={(0,0.4)}]
            \node at (0.5,1.5) {$\emptyset$};
      \end{tikzpicture} &
      \begin{tikzpicture}[scale=0.4,baseline={(0,0.4)}]
            \node at (0.5,1.5) {$\emptyset$};
      \end{tikzpicture} &
      \begin{tikzpicture}[scale=0.4,baseline={(0,0.4)}]
            \node at (0.5,1.5) {$\emptyset$};
      \end{tikzpicture} &	
	\begin{tikzpicture}[scale=0.4,baseline={(0,0)}]
		\draw (0,0) grid (1,1);
            \node at (0.5,0.5) {\tableaufontsize$5$};
       \end{tikzpicture} &
       \begin{tikzpicture}[scale=0.4,baseline={(0,0)}]
		\draw (0,0) grid (1,1);
            \node at (0.5,0.5) {\tableaufontsize$5$};
       \end{tikzpicture} &
       \begin{tikzpicture}[scale=0.4,baseline={(0,0)}]
		\draw (0,0) grid (1,1);
            \node at (0.5,0.5) {\tableaufontsize$5$};
       \end{tikzpicture} &
       \begin{tikzpicture}[scale=0.4,baseline={(0,0)}]
		\draw (0,0) grid (1,1);
            \node at (0.5,0.5) {\tableaufontsize$5$};
       \end{tikzpicture} &	
        \begin{tikzpicture}[scale=0.4,baseline={(0,0)}]
		\draw (0,0) grid (2,1);
            \node at (0.5,0.5) {\tableaufontsize$5$};
            \node at (1.5, 0.5) {\tableaufontsize$6$};
       \end{tikzpicture} &	
        \begin{tikzpicture}[scale=0.4,baseline={(0,0)}]
		\draw (0,0) grid (2,1);
            \node at (0.5,0.5) {\tableaufontsize$5$};
            \node at (1.5, 0.5) {\tableaufontsize$6$};
       \end{tikzpicture} &	
         \begin{tikzpicture}[scale=0.4,baseline={(0,0)}]
		\draw (0,0) grid (2,1);
            \node at (0.5,0.5) {\tableaufontsize$5$};
            \node at (1.5, 0.5) {\tableaufontsize$6$};
       \end{tikzpicture} &		
	\begin{tikzpicture}[scale=0.4,baseline={(0,0.4)}]
		\draw (0,1) grid (1,2);
            \node at (0.5,1.5) {\tableaufontsize$6$};
	\end{tikzpicture} &
	\begin{tikzpicture}[scale=0.4,baseline={(0,0.4)}]
		\draw (0,1) grid (1,2);
            \node at (0.5,1.5) {\tableaufontsize$6$};
	\end{tikzpicture}&
	\begin{tikzpicture}[scale=0.4,baseline={(0,0.4)}]
            \node at (0.5,1.5) {$\emptyset$};
      \end{tikzpicture} 
\\
Rule $H_2$ &
   \begin{tikzpicture}
      \draw (1.2,0) node[pnt]{};
   \end{tikzpicture}
   &
   \begin{tikzpicture}
      \draw (1.2,0) node[pnt]{};
   \end{tikzpicture}
   &
   \begin{tikzpicture}
      \draw[bend left] (1,0) node[pnt]{} to +(0.3,0.2);
   \end{tikzpicture}
   & 
   \begin{tikzpicture}
      \draw (1.2,0) node[pnt]{};
   \end{tikzpicture}
   &
  \begin{tikzpicture}
      \draw (1.2,0) node[pnt]{};
   \end{tikzpicture}
   &
   \begin{tikzpicture}
      \draw (1.2,0) node[pnt]{};
   \end{tikzpicture}
   & 
     \begin{tikzpicture}
      \draw[bend left] (1,0) node[pnt]{} to +(0.3,0.2);
   \end{tikzpicture}
   & 
    \begin{tikzpicture}
      \draw (1.2,0) node[pnt]{};
   \end{tikzpicture}
   &  
     \begin{tikzpicture}
      \draw (1.3,0) node[pnt]{};
   \end{tikzpicture}
   &
   \begin{tikzpicture}
      \draw[bend right] (1.5,0) node[pnt]{} to (1.2,0.2);
   \end{tikzpicture}
   &
   \begin{tikzpicture}
      \draw (1.3,0) node[pnt]{};
   \end{tikzpicture}
   &
   \begin{tikzpicture}
      \draw[bend right] (1.5,0) node[pnt]{} to (1.2,0.2);
   \end{tikzpicture}
\end{tabular}
\\
	\begin{tikzpicture}[scale=1.5]
		 \foreach \i in {1,...,6}
       	 \node[pnt,label=below:$\i$] at (\i,0)(\i) {};
       \draw[bend left=45](1)  to node[pos=0.2,above]{$\scriptstyle1$}  (4);
      \draw[bend left=30](4) to node[pos=0.7, above]{$\scriptstyle2$}  (6);
      \draw[bend left=45](2) to node[above]{$\scriptstyle2$} (5);
      \draw[loop above](3)  to node[pos=0.1,right]{$\scriptstyle1$}  ();
      \draw[bend right=30](2) to node[above]{$\scriptstyle2$}  (5);
      \draw[bend right=30](1) to node[below]{$\scriptstyle2$} (6);
   \end{tikzpicture}
\\
\begin{tabular}{*{13}{c}}
Rule $V_2$ &
   \begin{tikzpicture}[baseline={(1,-0.2)}]
      \draw (1,0) node[pnt]{};
   \end{tikzpicture}
   & 
   \begin{tikzpicture}[baseline={(1,-0.2)}]
      \draw[bend right] (1,0) node[pnt]{} to (1.3,-0.2);
   \end{tikzpicture}
   & 
    \begin{tikzpicture}[baseline={(1,-0.2)}]
      \draw (1,0) node[pnt]{};
   \end{tikzpicture}
   & 
   \begin{tikzpicture}[baseline={(1,-0.2)}]
      \draw[bend right] (1,0) node[pnt]{} to (1.3,-0.2);
   \end{tikzpicture}
   & 
    \begin{tikzpicture}[baseline={(1,-0.2)}]
      \draw (1,0) node[pnt]{};
   \end{tikzpicture}
   & 
    \begin{tikzpicture}[baseline={(1,-0.2)}]
      \draw (1,0) node[pnt]{};
   \end{tikzpicture}
   &   
   \begin{tikzpicture}[baseline={(1,-0.2)}]
      \draw (1,0) node[pnt]{};
   \end{tikzpicture}
   & 
    \begin{tikzpicture}[baseline={(1,-0.2)}]
      \draw (1,0) node[pnt]{};
   \end{tikzpicture}
   & 
   \begin{tikzpicture}[baseline={(1,-0.2)}]
      \draw[bend left] (1.3,0) node[pnt]{} to (1,-0.2);
   \end{tikzpicture}
   &
   \begin{tikzpicture}[baseline={(1,-0.2)}]
      \draw (1.5,0) node[pnt]{};
   \end{tikzpicture}
   &
   \begin{tikzpicture}[baseline={(1,-0.2)}]
      \draw[bend left] (1.3,0) node[pnt]{} to (1,-0.2);
   \end{tikzpicture}
   &
   \begin{tikzpicture}[baseline={(1,-0.2)}]
      \draw (1.5,0) node[pnt]{};
   \end{tikzpicture}
       \\[1.5ex]
    $\mu_2^0$ & $\mu_2^1$ & $\mu_2^2$ & $\mu_2^3$ & $\mu_2^4$ & $\mu_2^5$  & $\mu_2^6$ & $\mu_2^7$   & $\mu_2^8$  & $\mu_2^9$  & $\mu_2^{10}$   & $\mu_2^{11}$  & $\mu_2^{12}$
  \\
	\begin{tikzpicture}[scale=0.4,baseline={(0,0.4)}]
            \node at (0.5,1.5) {$\emptyset$};
      \end{tikzpicture} &
	\begin{tikzpicture}[scale=0.4,baseline={(0,0.4)}]
            \node at (0.5,1.5) {$\emptyset$};
      \end{tikzpicture} &
	\begin{tikzpicture}[scale=0.4,baseline={(0,0.4)}]
		\draw (0,1) grid (1,2);
		\node at (0.5,1.5) {\tableaufontsize$6$};
	\end{tikzpicture} &
	\begin{tikzpicture}[scale=0.4,baseline={(0,0.4)}]
		\draw (0,1) grid (1,2);
		\node at (0.5,1.5) {\tableaufontsize$6$};
	\end{tikzpicture} &
	\begin{tikzpicture}[scale=0.4,baseline={(0,0.4)}]
		\draw (0,0) grid (1,2);
		\node at (0.5,1.5) {\tableaufontsize$5$};
		\node at (0.5,0.5) {\tableaufontsize$6$};
	\end{tikzpicture} &
	\begin{tikzpicture}[scale=0.4,baseline={(0,0.4)}]
		\draw (0,0) grid (1,2);
		\node at (0.5,1.5) {\tableaufontsize$5$};
		\node at (0.5,0.5) {\tableaufontsize$6$};
	\end{tikzpicture} &
	\begin{tikzpicture}[scale=0.4,baseline={(0,0.4)}]
		\draw (0,0) grid (1,2);
		\node at (0.5,1.5) {\tableaufontsize$5$};
		\node at (0.5,0.5) {\tableaufontsize$6$};
	\end{tikzpicture} &
	\begin{tikzpicture}[scale=0.4,baseline={(0,0.4)}]
		\draw (0,0) grid (1,2);
		\node at (0.5,1.5) {\tableaufontsize$5$};
		\node at (0.5,0.5) {\tableaufontsize$6$};
	\end{tikzpicture} &
	\begin{tikzpicture}[scale=0.4,baseline={(0,0.4)}]
		\draw (0,0) grid (1,2);
		\node at (0.5,1.5) {\tableaufontsize$5$};
		\node at (0.5,0.5) {\tableaufontsize$6$};
	\end{tikzpicture} &
	\begin{tikzpicture}[scale=0.4,baseline={(0,0.4)}]
		\draw (0,1) grid (1,2);
		\node at (0.5,1.5) {\tableaufontsize$6$};
	\end{tikzpicture} &
	\begin{tikzpicture}[scale=0.4,baseline={(0,0.4)}]
		\draw (0,1) grid (1,2);
		\node at (0.5,1.5) {\tableaufontsize$6$};
	\end{tikzpicture} &
	\begin{tikzpicture}[scale=0.4,baseline={(0,0.4)}]
            \node at (0.5,1.5) {$\emptyset$};
      \end{tikzpicture} &
	\begin{tikzpicture}[scale=0.4,baseline={(0,0.4)}]
            \node at (0.5,1.5) {$\emptyset$};
      \end{tikzpicture} 
\end{tabular}
\end{center}
\end{example}
The result of transposing every tableau in each sequence $\lambda_1$, $\lambda_2$, and $\mu_2$, and filling the tableau from the right is the following $2$-coloured permutation in Figure~\ref{fig:inperm}.

\section{Enumeration of \texorpdfstring{$r$}{r}-coloured permutations}
\label{sec:enum}

Before we enumerate $r$-coloured permutations, a quick overview of Marberg's approach for the enumeration of coloured set partitions helps set the stage for a new interpretation.

\subsection{Another interpretation of 
   \texorpdfstring{$\mathcal{G}_{j,k,r}$}{Gjkr} for set partitions}

Marberg viewed $r$ sequences of vacillating tableaux, one for each colour, as $r \times (k-1)$ matrices $A = [A_{i,l}]$  encoding $\lambda^l_i$ in a vacillating tableau sequence $T$ for colour $i$. If the set partition is $j$-noncrossing and $k$-nonnesting, then this tableau has a maximum of $j-1$ columns and $k-1$ rows. For colour $i$, the $i$th row of matrix $A$ just lists parts of $\lambda^l$, thus at most $k-1$ non-zero parts. The multigraph $\mathcal{G}_{j,k,r}$ is drawn using all such allowable $A$'s as vertices, and edges and loops connecting vertices corresponding to adding a box, deleting a box, or doing nothing in the construction of vacillating tableaux so that the resulting sequence contains only tableaux of at most $j-1$ columns and $k-1$ rows. Once completed, the multigraph  $\mathcal{G}_{j,k,r}$ gives rise to an adjacency matrix. To find the number $\NCN_{j,k}(n,r)$ which is also the number of  ($n-1$)-step walks on  $\mathcal{G}_{j,k,r}$ from the zero matrix to itself, the method of transfer matrix gives a quotient of two polynomials (determinants actually), thus concluding that the ordinary generating function $\sum_{n \ge 0}\NCN_{j,k}(n+1,r) x^n$ is rational.

\subsection{Examples of 
\texorpdfstring{$\mathcal{G}_{2, 2, 1}$}{G221} 
and 
\texorpdfstring{$\mathcal{G}_{2, 2, 2}$}{G222} 
for set partitions}

To illustrate the construction of $\mathcal{G}_{j,k,r}$, we first reconstruct Marberg's $\mathcal{G}_{2,2,1}$ and $\mathcal{G}_{2,2,2}$ by naming each vertex and edge as it becomes necessary.

The arc annotated diagram of a set partition on $[n]$ has $n-1$ consecutive gaps, i.\,e.\ between each pair of adjacent points. Let the set of non-crossing, non-nesting, uncoloured set partitions on $[n]$ be denoted by $\mathcal{P}_{2,2,1}(n)$. For each $P \in \mathcal{P}_{2,2,1}(n)$, a snap shot of each gap belongs to one of the first four types in Table~\ref{tab:five} where the matching steps in $\mathcal{G}_{2,2,1}$ are also given. Since $r=1$, only two vertices exist in $\mathcal{G}_{2,2,1}$: $v_0$, the initial state for no opener, and $v_1$, for one opener. No other vertices accounting for other states are present because any state $v_i$ where $i \ge 2$ would mean two or more openers which will form at least a $2$-nesting or $2$-crossing when closed. Incident at $v_0$ are three types of edges: two loops, 
\begin{tikzpicture}[baseline={(0,0)}]
      \draw (0,0) node[pnt]{} to[loop above] node[right]{$\scriptstyle\times$} ();
\end{tikzpicture}
  for no arc in the consecutive gap, and 
 \begin{tikzpicture}[baseline={(0,0)}]
      \draw (0,0) node[pnt]{} to[loop above] node[right]{$\scriptstyle1$} ();
 \end{tikzpicture}
 for a distance $1$-arc both of which do not change the number of openers present as the set partition is scanned from the left to the right; the last type is a directed edge from $v_0$ to $v_1$ to indicate that an opener is present in the consecutive gap. Once at $v_1$, only the loop,
 \begin{tikzpicture}[baseline={(0,0)}]
      \draw (0,0) node[pnt]{} to[loop above] node[right]{$\scriptstyle\times$} ();
 \end{tikzpicture},
 is allowed because a $1$-arc
 \begin{tikzpicture}[baseline={(0,0)}]
      \draw (0,0) node[pnt]{} to[loop above] node[right]{$\scriptstyle1$} ();
   \end{tikzpicture}
will create a $2$-nesting in $P$ with the existing opener. A directed edge from $v_1$ to $v_0$ means that an opener is closed. To simplify drawing, an edge without arrows is bidirectional. The result is shown in Figure~\ref{fig:setpartG221}.

The adjacency matrix of Figure~\ref{fig:setpartG221} is 
\[
\left[\begin{matrix} 2, 1\\ 1, 1\end{matrix} \right]
\]
which gives 
\[
\sum_{n \ge 0} \NCN_{2,2}(n+1, 1) x^n = \frac{1-x}{1-3x+x^2}
\]
and expands to
\begin{quote}\raggedright
\(
1+2x+5x^2+13x^3+34x^4+89x^5+233x^6+610x^7+1597x^8+4181x^9
+10946x^{10}+28657x^{11}+75025x^{12}+196418x^{13}+514229x^{14}
+1346269x^{15}+3524578x^{16}+9227465x^{17}+24157817x^{18}
+63245986x^{19}+O(x^{20}).
\)
\end{quote}
As noted by \citet{Mar12}, these are every second Fibonacci numbers.

To construct  $\mathcal{G}_{2,2,2}$, we require four vertices: still $v_0$ as the initial state for no  opener, but also two states indicating one $r$-coloured ($r \in [2]$) opener, $v_{1_1}$ and $v_{1_2}$. Since two arcs of different colours do not create a crossing or nesting, one more state is needed, $v_{2_{12}}$, for two openers, one of each colour. As in  $\mathcal{G}_{2,2,1}$, the loops and edges are placed according to what is allowed in $P$, but a new edge between $v_{1_1}$ and $v_{1_2}$ is added in the last row of Table~\ref{tab:five} for the closing of one colour on point $m$ while an opener is present at point $m-1$ in $P$. The result is shown in Figure~\ref{fig:setpartG222} with its associated adjacency matrix
\[
\left[
\begin{matrix}
	3,1,1,0\\
	1,2,1,1\\
	1,1,2,1\\
	0,1,1,1
\end{matrix}
\right]
\]
and generating function
\[
\sum_{n \ge 0} \NCN_{2,2}(n+1, 2) x^n = \frac{1-4x + x^2}{1-7x+11x^2 - x^3}
\]
which expands to
\begin{quote}\raggedright
\(
1+3x+11x^2+45x^3+197x^4+895x^5+4143x^6+19353x^7
+90793x^8+426811x^9+2008307x^{10}+9454021x^{11}+44513581x^{12}
+209609143x^{13}+987068631x^{14}+4648293425x^{15}+21889908177x^{16}
+103085198195x^{17}+485455690843x^{18}+2286142563933x^{19}+O(x^{20})
\)
\end{quote}
Once automated, we obtain new series easily, A225029 and A225030 in \citep{oeis} are for $3$ and $4$-coloured set partitions.
\begin{multline*}
\sum_{n \ge 0} \NCN_{2,2}(n, 3) x^n =
                                 \frac{1-10x+22x^2-x^3}{1-14x +59 x^2-74x^3+x^4}\\
=1+4x+19x^2+103x^3+616x^4+3949x^5+26545x^6 + 184120x^7+ O(x^8).
\end{multline*}
\begin{multline*}
\sum_{n \ge 0} \NCN_{2,2}(n, 4) x^n = 
          \frac{1-20x+122x^2-224x^3+x^4}{1-25x +218 x^2-782x^3+973x^4-x^5}\\
=1+5x+29x^2+193x^3+1441x^4+11765x^5 +102701x^6+ 941857x^7 + O(x^8).
\end{multline*}
Using an average personal computer, Maple15 can generate up to $7$ colours. See A225031--A225033 in \citep{oeis}. The next case with a matrix size of $256 \times 256$, computation would take too long to find the determinants.

\begin{figure}
\begin{minipage}{0.45\linewidth}
\figurefontsize\centering
\begin{tikzpicture}
	\node[state] (O) at (0,0) {$v_0$};
	\node[state] (I) at (2,0) {$v_1$};
	\draw[loop  above left] (O) to node[above]{$\scriptstyle\times$}();
	\draw[loop  below left] (O) to node[below]{$\scriptstyle1$}();	
	\draw[loop right] (I) to node[right]{$\scriptstyle\times$}();
	\draw (O) to node[above]{$\scriptstyle1$}(I);
\end{tikzpicture}
\caption{An uncoloured set partition graph, $\mathcal{G}_{2, 2, 1}$.}
\label{fig:setpartG221}
\end{minipage}
\qquad
\begin{minipage}{0.45\linewidth}
\figurefontsize\centering
\begin{tikzpicture}[state/.append style={minimum size=8mm},scale=0.8]
	\node[state] (O) at (-1,0) {$v_0$};
	\node[state] (A) at (1,1) {$v_{1_1}$};
	\node[state] (B) at (1, -1) {$v_{1_2}$};
	\node[state] (C) at (3, 0) {$v_{2_{12}}$};
	\draw[loop  above] (O) to node[above]{$\scriptstyle\times$} ();
	\draw[loop  below] (O) to node[below]{$\scriptstyle1$}();
	\draw[loop  left] (O) to node[left]{$\scriptstyle2$}();	
	\draw[loop above left] (A) to node[left]{$\scriptstyle\times$} ();
	\draw[loop above right] (A) to node[right]{$\scriptstyle2$}();
	\draw[loop below left] (B) to node[left]{$\scriptstyle\times$}();
	\draw[loop below right] (B) to node[below]{$\scriptstyle1$}();
	\draw[loop right] (C) to node[right]{$\scriptstyle\times$}();
	\draw (O) to node[above]{$\scriptstyle1$} (A);
	\draw (O) to node[below]{$\scriptstyle2$} (B);
	\draw (A) to node[right]{$\scriptstyle12$} (B);
	\draw (A) to node[right]{$\scriptstyle2$} (C);
	\draw (B) to node[right]{$\scriptstyle1$} (C);	
\end{tikzpicture}
\caption{A $2$-coloured set partition graph, $\mathcal{G}_{2, 2, 2}$.}
\label{fig:setpartG222}
\end{minipage}
\end{figure}

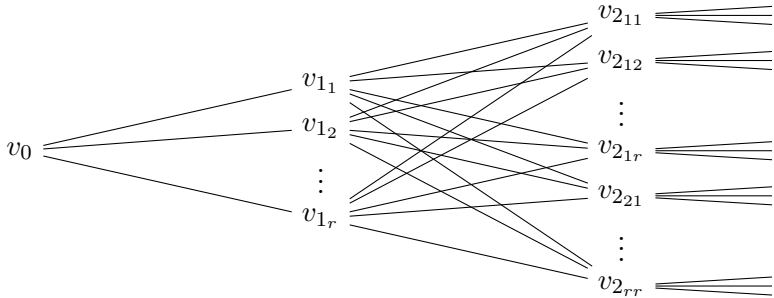
\begin{figure}
\centering\figurefontsize
\begin{tikzpicture}[xscale=4,yscale=0.6]
   \node (0) at (0,0) {$v_0$};
   \node (11) at (1,1.5) {$v_{1_1}$};
   \node (12) at (1,0.5) {$v_{1_2}$};
   \node at (1,-0.5) {$\vdots$};
   \node (1r) at (1,-1.5) {$v_{1_r}$};
   \node (211) at (2,3) {$v_{2_{11}}$};
   \node (212) at (2,2) {$v_{2_{12}}$};
   \node at (2,1) {$\vdots$};
   \node (21r) at (2,0) {$v_{2_{1r}}$};
   \node (221) at (2,-1) {$v_{2_{21}}$};
   \node at (2,-2) {$\vdots$};
   \node (2rr) at (2,-3) {$v_{2_{rr}}$};
   \foreach \i in {1,2,r} {
      \draw (0) -- (1\i);
      \foreach \j in {11,12,1r,21,rr}
         \draw (1\i) -- (2\j);
   }
   \foreach \j in {11,12,1r,21,rr} {
      \draw (2\j) -- +(0.5,0.2);
      \draw (2\j) -- +(0.5,0);
      \draw (2\j) -- +(0.5,-0.2);
   }
\end{tikzpicture}
\caption{The line-up for states of the same number of openers}
\label{fig:vertices}
\end{figure}

In general, we obtain  $\mathcal{G}_{j,k,r}$ directly through labelling the edges and vertices of  $\mathcal{G}_{j,k,r}$ similar to generating such set partitions through the method of generating trees except that each vertex $v_i$ (considered as a state) in $\mathcal{G}_{j,k,r}$ indicates that $i$ openers are pending to close. When drawn from the left to the right where all vertices of the same first subscript line up vertically, we get edges either between $v_i$ and $v_{i+1}$ for each $i \ge 0$ for openers or closers  as in Figure~\ref{fig:vertices}, or between vertices of the same first subscript for the presence of both (drawn as vertical edges, not shown in  Figure~\ref{fig:vertices}). Care needs to be taken when many arcs of the same colour are open because the order in which they are closed relates to how crossing and nesting are formed.

\begin{table}[htb]
\figurefontsize\centering
\begin{tabular}{cccc}
Domain & Set partition & Types of arcs & 
\multicolumn{1}{p{3cm}}{Types of steps in $\mathcal{G}_{2,2,1}$ and  $\mathcal{G}_{2,2,2}$}
\\ \hline
   $m \ge 2$ 
   &
   \begin{tikzpicture}
      \path (0,0) node[pnt]{} node[below]{$m-1$}
               (1,0) node[pnt]{} node[below]{$m\phantom{1}$};
   \end{tikzpicture}
   & no arc & 
   \begin{tikzpicture}[baseline={(0,0)}]
      \draw (0,0) node[pnt]{} to[loop above] node[right]{$\scriptstyle\times$} ();
   \end{tikzpicture}
\\
   $m \ge 2$ 
   &
   \begin{tikzpicture}
      \draw[bend left] (0,0) node[pnt]{} node[below]{$m-1$}
            to node[above]{$\scriptstyle r$}
            (1,0) node[pnt]{} node[below]{$m\phantom{1}$};
   \end{tikzpicture}
   & a $1$-arc coloured $r$  &  
   \begin{tikzpicture}
      \draw (0,0) node[pnt]{} to[loop above] node[right]{$\scriptstyle r$} ();
   \end{tikzpicture},
   $r \in [2]$
\\
   $m \ge 2$ 
   &  
   \begin{tikzpicture}
      \draw[bend left,uni] (0,0) node[pnt]{} node[below]{$m-1$}
            to node[above]{$\scriptstyle r$}
            +(0.7,0.3);
      \draw (1,0) node[pnt]{} node[below]{$m\phantom{1}$};
   \end{tikzpicture}
   & an opener & 
   \begin{tikzpicture}
      \draw[uni] (0,0) node[pnt]{} node[below]{$v_i$}
         -- node[above]{$\scriptstyle r$}
         (1,0) node[pnt]{} node[below]{$v_{i+1}$};
   \end{tikzpicture}
\\
   $m \ge 3$ 
   &  
   \begin{tikzpicture}
      \draw[bend right,unir] (1,0) node[pnt]{} node[below]{$m\phantom{1}$}
            to node[above]{$\scriptstyle r$}
            +(-0.7,0.3);
      \draw (0,0) node[pnt]{} node[below]{$m-1$};
   \end{tikzpicture}
   & a closer  & 
   \begin{tikzpicture}
      \draw[uni] (1,0) node[pnt]{} node[below]{$v_{i+1}$}
         -- node[below]{$\scriptstyle r$}
         (0,0) node[pnt]{} node[below]{$v_{i}$};
   \end{tikzpicture}
\\
   $m \ge 3$ 
   &  
   \begin{tikzpicture}
      \draw[bend left,uni] (0,0) node[pnt]{} node[below]{$m-1$}
            to node[above,near end]{$\scriptstyle 2$}
            +(1,0.3);
      \draw[bend right,unir] (1,0) node[pnt]{} node[below]{$m\phantom{1}$}
            to node[above,near end]{$\scriptstyle 1$}
            +(-1,0.3);
   \end{tikzpicture}
   & a closer and an opener 
   &  
   \begin{tikzpicture}
      \draw[uni] (0,1) node[pnt]{} node[right]{$v_{i_1}$}
         -- node[right] {$\scriptstyle 12$}
         (0,0) node[pnt] {} node[right]{$v_{i_2}$};
   \end{tikzpicture}
\end{tabular}
\caption[Set partition and its multigraph]
{Five situations between point $m-1$ and point $m$ for set partitions and the matching steps in $\mathcal{G}$.}
\label{tab:five}
\end{table}

\subsection{Multigraphs, 
\texorpdfstring{$\mathcal{G}_{2, 2, 1}$}{G221} 
and 
\texorpdfstring{$\mathcal{G}_{2, 2, 2}$}{G222}
for permutations}

Instead of translating consecutive gaps  from set partitions into steps in the multigraph $\mathcal{G}$,  we examine each vertex in the arc diagram of a coloured permutation and assign each type of vertex to a step in $\mathcal{G}$. As for set partitions, we first construct the multigraph $\mathcal{G}_{2,2,1}$ for non-crossing, non-nesting, uncoloured permutations. Let us denote the set of all such permutations on $[n]$ by $\mathcal{S}_{2,2,1}(n)$. If $S \in \mathcal{S}_{2,2,1}(n)$, then a vertex is either
 a fixed point (
 \begin{tikzpicture}
      \draw[loop above] (0,0) node[pnt]{} to (0.0);
   \end{tikzpicture}
   ) , 
an opener (
\begin{tikzpicture}
      \draw[bend left] (1.2,0) node[pnt]{} to +(0.3,0.2);
      \draw[bend right] (1.2, 0) to (1.5, -0.2);
\end{tikzpicture}
),
 a closer (
  \begin{tikzpicture}
      \draw[bend right] (1.5,0) node[pnt]{} to (1.2, 0.2);
      \draw[bend left] (1.5, 0) to (1.2, -0.2);
   \end{tikzpicture}
 ), 
or   a lower transitory(
   \begin{tikzpicture}
      \draw[bend right] (0,0) node[pnt]{} to (0.3,-0.2);
      \draw[bend left] (0,0) to (-0.3, -0.2);
   \end{tikzpicture}
   ). We can't have an upper transitory which contributes to a $2$-(enhanced) crossing.

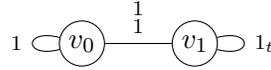
\begin{figure}
\figurefontsize\centering
\begin{tikzpicture}
	\node[state] (O) at (0,0) {$v_0$};
	\node[state] (I) at (1.5,0) {$v_1$};
	\draw[loop  left] (O) to node[left]{$\scriptstyle1$}();	
	\draw[loop right] (I) to node[right]{$\scriptstyle1_t$}();
	\draw (O) to node[above]{$\begin{smallmatrix}1\\1\end{smallmatrix}$}(I);
\end{tikzpicture}
\caption{An uncoloured permutation graph, $\mathcal{G}_{2, 2, 1}$.}
\label{fig:permG221}
\end{figure}

In Figure~\ref{fig:permG221}, $v_0$ still indicates the initial state with $0$ opener; $v_1$ indicates the state with $1$ opener. The loop labelled $1$ is the step taken when a fixed point coloured $1$ is encountered in the permutation scanned from the left. The loop labelled $1_t$ is the presence of a lower transitory with coloured $1$ arcs on both sides; this is possible only when an opener coloured $1$ is present, thus at $v_1$. Note that a lower transitory does not alter the state. The directed edge $(v_0, v_1)$ indicates the presence of an opener, and the edge traversed in reverse indicates that of a closer. An edge drawn without arrows still means a bidirectional edge.

The construction of $\mathcal{G}_{2,2,2}$ involves more types of vertices and edges which we summarize in Table~\ref{tab:perm}. Each state with one opener has the colours of the openers as subscripts. When a state has two openers, both colours are used, thus only one such vertex in $v_2$. The method of transfer matrix gives the following generating function. Here $x$ marks the size of the permutation.
\[
\sum_{n \ge 0} \NCN_{2,2}(n, 2) x^n =  \frac{1-6x+4x^2}{(1-2x)(1-6x)}
\]
which generates
\begin{quote}\raggedright
\(
1+2x+8x^2+40x^3+224x^4+1312x^5+7808x^6+46720x^7
+280064x^8+1679872x^9+10078208x^{10}
+60467200x^{11}+362799104x^{12}
+2176786432x^{13}+13060702208x^{14}+O(x^{15}).
\)
\end{quote}
The first group of non-trivial noncrossing and nonnesting, $2$-coloured permutations is on $[4]$. For each of the $24$ permutations on $[4]$, $8$ can be coloured in $4$ ways each; $8$, in $8$ ways each; and $8$, in $16$ ways each; thus, $8 \times (4 + 8 + 16) = 224$. This series, A092807 in \citep{oeis}, counts (with interpolated zeros) the number of closed walks of length $n$ at a vertex of the edge-vertex incidence graph of  $K_4$, the complete graph on $4$ vertices associated with the edges of $K_4$.

The next two series, A224992 and A224993 in \citep{oeis}, however, are new, namely,  $3$ and $4$-coloured noncrossing, nonnesting permutations. 
For $5$ colours, the matrix size, $252 \times 252$, hinders fast computation of determinants.
\begin{multline*}
\sum_{n \ge 0} \NCN_{2,2}(n, 3) x^n =
                                                \frac{1-17x+66x^2-36x^3}{(1-2x)(1-6x)(1-12x)}\\
=1+3x+18x^2+144x^3+1368x^4+14400x^5 + 160992x^6 + 1861632x^7+O(x^8),
\end{multline*}
\begin{multline*}
\sum_{n \ge 0} \NCN_{2,2}(n, 4) x^n = 
                       \frac{1-36x+380x^2-1200x^3+576x^4}{(1-2x)(1-6x)(1-12x)(1-20x)}\\
=1+4x+32x^2+352x^3+4736x^4+72832x^5+ 1226240x^6 + 21948928x^7 + O(x^8).
\end{multline*}

\begin{table}[htb]
\figurefontsize\centering
\setlength{\tabcolsep}{3pt}
\begin{tabular}{cccc}
Domain & Permutation & Types of arcs & Types of steps in  $\mathcal{G}_{2,2,2}$
\\
            & Vertex\\ \hline
  all vertices
   &
   \begin{tikzpicture}
      \draw[loop above] (0,0) node[pnt]{} to node[above]{$\scriptstyle l$}(0.0);
   \end{tikzpicture}
   & a fixed point & 
   \begin{tikzpicture}[baseline={(0,0)}]
      \draw (0,0) node[pnt]{} to[loop above] node[right]{$\scriptstyle l$} ();
   \end{tikzpicture},
   $l \in [2]$
\\
  all  except the last
   &
   \begin{tikzpicture}
      \draw[bend left] (1.2,0) node[pnt]{} 	
      					to node[above]{$\scriptstyle r$}
      					+(0.3,0.2);
      \draw[bend right] (1.2, 0) 
      					to node[below]{$\scriptstyle s$}(1.5, -0.2);
\end{tikzpicture}
   & an opener   &  
    \begin{tikzpicture}
      \draw[uni] (0,0) node[pnt]{} node[below]{$v_i$}
         -- node[above]{$\begin{smallmatrix}r\\s\end{smallmatrix}$}
         (1,0) node[pnt]{} node[below]{$v_{i+1}$};
   \end{tikzpicture},
   $r, s \in [2]$
\\
   
  all except the first 
   &  
   \begin{tikzpicture}
      \draw[bend right] (1.5,0) node[pnt]{} 
      					to node[above]{$\scriptstyle r$}
      					 (1.2, 0.2);
      \draw[bend left] (1.5, 0) 
      					to node[below]{$\scriptstyle s$}(1.2, -0.2);
   \end{tikzpicture}
   & a closer  & 
   \begin{tikzpicture}
      \draw[uni] (1,0) node[pnt]{} node[below]{$v_{i+1}$}
         -- node[above]{$\begin{smallmatrix}r\\s\end{smallmatrix}$}
         (0,0) node[pnt]{} node[below]{$v_{i}$};
   \end{tikzpicture},
   $r, s \in [2]$
\\
   no first, no last 
   &  
  \begin{tikzpicture}
      \draw[bend right] (0,0) node[pnt]{} 
      					to node[below]{$\scriptstyle r$}
      					 (0.3,-0.2);
      \draw[bend left] (0,0) 
      					to node[below]{$\scriptstyle r$}
      					(-0.3, -0.2);
   \end{tikzpicture}
   & a lower transitory
   &  
    \begin{tikzpicture}[baseline={(0,0)}]
      \draw (0,0) node[pnt]{} to[loop above] node[right]{$\scriptstyle r_t$} ();
   \end{tikzpicture},
   $r \in [2]$
\\
no first, no last
&
\begin{tikzpicture}
      \draw[bend left] (0,0) node[pnt]{} 
      					to node[above]{$\scriptstyle s$}+(0.3,0.2);
      \draw[bend right] (0,0) 
      					to node[above]{$\scriptstyle r$}(-0.3, 0.2);
   \end{tikzpicture}
 &an  upper transitory
 &
  \begin{tikzpicture}
      \draw[uni] (0,1) node[pnt]{} node[right]{$v_{i_r}$}
         -- node[right] {$\scriptstyle \underline{rs}$}
         (0,0) node[pnt] {} node[right]{$v_{i_s}$};
   \end{tikzpicture},
   $r, s \in [2]$
\\
   no first, no last 
   &  
  \begin{tikzpicture}
      \draw[bend right] (0,0) node[pnt]{} 
      					to node[below]{$\scriptstyle r$}
      					 (0.3,-0.2);
      \draw[bend left] (0,0) 
      					to node[below]{$\scriptstyle s$}
      					(-0.3, -0.2);
   \end{tikzpicture}
   & a lower transitory
   &  
    \begin{tikzpicture}
      \draw[unir] (0,1) node[pnt]{} node[right]{$v_{i_r}$}
         -- node[right] {$\scriptstyle \overline{sr}$}
         (0,0) node[pnt] {} node[right]{$v_{i_s}$};
   \end{tikzpicture},
   $r, s \in [2]$
\end{tabular}
\caption{Vertices in permutations and the matching steps in $\mathcal{G}_{2,2,2}$.
}
\label{tab:perm}
\end{table}

\begin{figure}
\label{fig:permG222}
\figurefontsize\centering
\begin{tikzpicture}
\node[state] (O) at (0,0) {$v_0$};
   \node[state] (A) at (4, 3) {$v_{1_{\begin{smallmatrix}1\\1\end{smallmatrix}}}$};
   \node[state] (B) at (4, 1) {$v_{1_{\begin{smallmatrix}1\\2\end{smallmatrix}}}$};
   \node[state] (C) at (4,-1) {$v_{1_{\begin{smallmatrix}2\\1\end{smallmatrix}}}$};
   \node[state] (D) at (4,-3) {$v_{1_{\begin{smallmatrix}2\\2\end{smallmatrix}}}$};
   \node[state] (E) at (8, 0) {$v_{2_{\begin{smallmatrix}1,2\\1,2\end{smallmatrix}}}$};
   \draw[loop above left]  (O) to node[above]{$\scriptstyle1$} ();
   \draw[loop below left]  (O) to node[below]{$\scriptstyle2$} ();
   \draw[loop above left]  (A) to node[above]{$\scriptstyle2$} ();
   \draw[loop above right] (A) to node[above]{$\scriptstyle1_t$} ();
   \draw[loop above left]  (B) to node[above right]{$\scriptstyle2$} ();
   \draw[loop above right] (B) to node[above]{$\scriptstyle2_t$} ();
   \draw[loop below right] (C) to node[below]{$\scriptstyle1_t$} ();
   \draw[loop below left]  (C) to node[below right]{$\scriptstyle1$} ();
   \draw[loop below left]  (D) to node[below]{$\scriptstyle1$} ();
   \draw[loop below right] (D) to node[below]{$\scriptstyle2_t$} ();
   \draw[loop above right]       (E) to node[right]{$\scriptstyle1_t$} ();
   \draw[loop below right]       (E) to node[right]{$\scriptstyle2_t$} ();
   \draw (O)
      to node[above]{$\begin{smallmatrix}1\\1\end{smallmatrix}$} (A)
      to node[right]{$\scriptstyle\overline{12}$} (B)
      to node[above]{$\begin{smallmatrix}1\\2\end{smallmatrix}$} (O)
      to node[below]{$\begin{smallmatrix}2\\1\end{smallmatrix}$} (C)
      to node[right]{$\scriptstyle\overline{12}$} (D)
      to node[below]{$\begin{smallmatrix}2\\2\end{smallmatrix}$} (O);
   \draw[bend right] (A)
      to node[near start,left]{$\scriptstyle\underline{12}$} (C);
   \draw[bend right] (B)
      to node[near end,left]{$\scriptstyle\underline{12}$} (D);
   \draw (A)
      to node[above]{$\begin{smallmatrix}2\\2\end{smallmatrix}$} (E)
      to node[above]{$\begin{smallmatrix}1\\2\end{smallmatrix}$} (C);
   \draw[postaction={decoration}] (E)
      to node[above]{$\begin{smallmatrix}2\\1\end{smallmatrix}$} (B);
   \draw (E)
      to node[below]{$\begin{smallmatrix}1\\1\end{smallmatrix}$} (D);
      to node[near start,above]{$\begin{smallmatrix}2\\2\end{smallmatrix}$} (A);
      to node[below]{$\begin{smallmatrix}1\\2\end{smallmatrix}$} (C);
\end{tikzpicture}
\caption[A $2$-coloured permutation multigraph]
   {A $2$-coloured permutation multigraph, $\mathcal{G}_{2, 2, 2}$
   }
\end{figure}
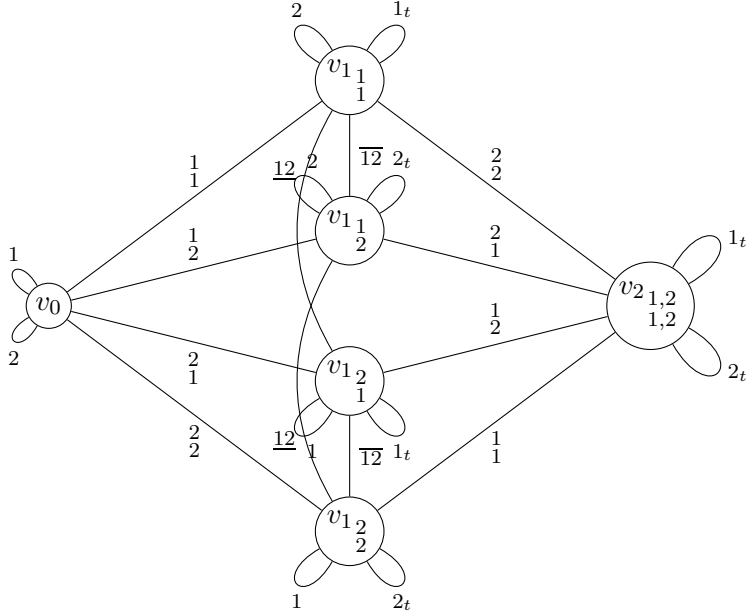

\subsection{Proof of Rationality through Multigraphs for \texorpdfstring{$r$}{r}-coloured permutations}

In general, drawing $\mathcal{G}_{j,k,r}$ for coloured permutations is a tedious task. As the $j$, $k$, and $r$ increase, types of edges and vertices increase. Not only does one need to track the order in which coloured arcs are closed, one also needs to create unidirectional edges which go to the right states. Regardless of the complexity of the multigraph, $\mathcal{G}_{j,k,r}$, only a finite number of vertices and edges are present because both crossing and nesting numbers are bounded for the set of $r$-coloured permutations. Furthermore, the number of such permutations on $[n]$ is the number of $n$-step paths from $v_0$ to $v_0$ in $\mathcal{G}_{j,k,r}$ because all openers must be closed. Using the method of transfer matrix then yields a rational function for the ordinary generating function, 
\(
\sum_{n \ge 0} \NCN_{j,k}(n,r) x^n.
\)

\section{Concluding Remarks}

\citet{ChG11} constructed oscillating $r$-rim hook tableaux for $r$-coloured complete matchings to establish symmetric joint distribution. This is equivalent to counting the number of walks on the Hasse diagram of the $r$-rim hook lattice. \citet{Mar12} noted isomorphism between his walks  on the Hasse diagram of the $r$-fold product of the Young lattice of integer partitions (produced from oscillating $r$-partite tableau) and those of \citet{ChG11}. Arc-coloured permutations have two types of $r$-partite tableau simultaneously accounting for both the upper and the lower arc diagrams. Does this correspond to a double $r$-rim hook lattice? 

When both nesting and crossing numbers are bounded, a finite multigraph can be constructed. This method of transfer matrix may be extended to the enumeration of set partitions of classical types as in the works of  \citet{RuSt10}, even their coloured counterparts. The challenge lies in finding the generating function when only one of the bounds is present. For instance, \citet{Mar12} showed that the ordinary generating function for noncrossing $2$-coloured set partitions is D-finite, but conjectured non-D-finite series for noncrossing $r$-coloured set partitions when $r \ge 3$.

\section{Acknowledgements}
The author would like to thank Marni Mishna for recommending Marberg's paper~\citep{Mar12},  Eric Marberg for clarifying the construction of the matrices which act as vertex labels of the multigraph $\mathcal{G}_{j,k,r}$, and two anonymous referees for their constructive comments which led to the automation of generating the rational series for coloured, noncrossing, nonnesting set partitions and permutations. Mogens Lemvig Hansen helped with figures and Maple code  for the automation.

\section{Implementation of automation}
Below is the Maple code that translates the graph, $\mathcal{G}_{2,2,c}$, of a $c$ coloured, noncrossing, nonnesting set partition to its adjacency matrix, $A$. 
\begin{enumerate}
\item Class$[k]$ indicates the vertices of the graph, each with its own state indicating the number of openers pending to close. 
\item The adjacency matrix $A$ has size $n\times n$ where $n$ is the number of vertices in the graph where the number of closed walks from $v_0$ is enumerated. 
\item The canonical ordering of the vertices of $\mathcal{G}$ begins with $v_0$ as the first vertex, followed by vertices of state $1$, ordered according to the subscript which indicates the colour of the open arc, and so on.
\item The diagonal entries of $A$ are the number of loops at each vertex of $\mathcal{G}_{2,2, c}$. 
\item The edges between vertices of consecutive states and vertices within the same state are described according to the size of the symmetric difference of the sets of openers, namely, subscripts of the vertices. These are translated to off diagonal entries of $A$.
\end{enumerate}

\begin{lstlisting}
setpartitionG22 := proc(c::posint, vlabel::name)
   local class, i, j, k, n, v, U, d, L, A;
   # The first vertex of class k is vertex[ class[k] ]
   class := array(0..c+1, [1]);
   for k from 0 to c do
      class[k+1] := class[k] + binomial(c, k);
   od;
   n := class[c+1] - 1;
   # v[i] is the label of vertex[i]
   v := array(1..n);
   i := 0;
   for k from 0 to c do
      for L in combinat['choose']( {$1..c}, k ) do
         i := i + 1;
         v[i] := L;
      od;
   od;
   if nargs >= 2 then vlabel := v fi;
   A := Matrix(n, n, 'shape'='symmetric');
   for k from 0 to c do
      for i from class[k] to class[k+1]-1 do
         # loops
         A[i,i] := 1 + c - k;
         # edges within class k
         for j from class[k] to i-1 do
            if nops( v[i] intersect v[j] ) = k-1 then
               A[i,j] := 1;
            fi;
         od;
         # edges from class k to class k+1
         if k = c then next fi;
         for j from class[k+1] to class[k+2]-1 do
            if  nops( v[j] minus v[i] ) = 1 then
               A[i,j] := 1;
            fi;
         od;
      od;
   od;
   A;
end:
\end{lstlisting}

Below is a similar Maple code for translating the graph, $\mathcal{G}_{2,2,c}$, of a $c$-coloured, noncrossing, nonnesting permutation to its adjacency matrix, $A$.
The construction parallel's that of coloured set partitions.

\begin{lstlisting}
permutationG22 := proc(c::posint, vlabel::name)
   local class, i, j, k, n, v, U, d, L, A;
   # The first vertex of class k is vertex[ class[k] ]
   class := array(0..c+1, [1,2]);
   for k from 1 to c do
      class[k+1] := class[k] + binomial(c, k)^2;
   od;
   n := class[c+1] - 1;
   # v[i] is the label of vertex[i]
   v := array(1..n);
   i := 0;
   for k from 0 to c do
      for U in combinat['choose']( {$1..c}, k ) do
         for L in combinat['choose']( {$1..c}, k ) do
            i := i + 1;
            v[i] := [ U, L ];
         od;
      od;
   od;
   if nargs >= 2 then vlabel := v fi;
   A := Matrix(n, n, 'shape'='symmetric');
   # loops
   for i from 1 to n do
      A[i,i] := c;
   od;
   for k from 0 to c do
      for i from class[k] to class[k+1]-1 do
         # edges within class k
         for j from class[k] to i-1 do
            if nops( v[i][1] intersect v[j][1] ) 
             + nops( v[i][2] intersect v[j][2] )
                = 2*k-1 then
               A[i,j] := 1;
            fi;
         od;
         # edges from class k to class k+1
         if k = c then next fi;
         for j from class[k+1] to class[k+2]-1 do
            if  nops( v[j][1] minus v[i][1] ) = 1
            and nops( v[j][2] minus v[i][2] ) = 1
            then
               A[i,j] := 1;
            fi;
         od;
      od;
   od;
   A;
end:
\end{lstlisting}

Once the adjacency matrix, $A$, is constructed, we use the linear algebra package to find the rational function via transfer matrix method. Note that by defining $A$ in the code to be permutationG22(c) instead of setpartitionG22(c), one will obtain the generating functions for $c$-coloured permutations. Depending on the version of Maple, earlier versions (around Maple 5) require that the determinant of the minor be spelled out as Determinant(Minor(B,1,1)). The code below runs in Maple 15.

\begin{lstlisting}
with(LinearAlgebra):
for c from 1 to 5 do 
   A := setpartitionG22(c);
   print(A);
   Id := IdentityMatrix( Dimension(A) ):
   B := Add(Id, A, 1, -lambda):
   r := Minor(B, 1,1) / Determinant(B):
   print(factor(r));
   print(series(r, lambda=0, 15));
od:
\end{lstlisting}
With an average personal computer running Windows 7, Maple 15 managed to find the rational functions up to $7$ colours for set partitions, then a $256 \times 256$ adjacency matrix took too long for $8$ colours to compute the ratio of determinants.

Coloured permutations are further restricted due to its complexity that already at $5$ colours, the adjacency matrix of size $252 \times 252$ was too large to compute its ratio of determinants.

\bibliographystyle{plainnat}

\end{document}